
\magnification=\magstep1
\overfullrule=0pt
     
\def\eqde{\,{\buildrel \rm def \over =}\,}
\def\al{\alpha} \def\la{\lambda} \def\be{\beta} \def\ga{\gamma}
\def\J{{\cal J}}    \def\z{{\cal Z}}  \def\G{\Gamma}
\def\si{\sigma}  
\def\ade{${\cal ADE}_7$-type invariant}   \def\va{\varphi}
\def\s{{\cal S}}    \def\k{\bar{k}}   \def\I{{\cal I}}
\def\r{\bar{r}}  \def\L{{\Lambda}} 
\def\E{{\cal E}}  \def\p{{\cal P}}  \def\D{{\cal D}} \def\F{{\cal F}}
\font\huge=cmr10 scaled \magstep2
\def\QED{\vrule height6pt width6pt depth0pt}
\font\smcap=cmcsc10     \font\little=cmr7

\catcode`\@=11
\font\tenmsa=msam10
\font\sevenmsa=msam7
\font\fivemsa=msam5
\font\tenmsb=msbm10
\font\sevenmsb=msbm7
\font\fivemsb=msbm5
\newfam\msafam
\newfam\msbfam
\textfont\msafam=\tenmsa  \scriptfont\msafam=\sevenmsa
  \scriptscriptfont\msafam=\fivemsa
\textfont\msbfam=\tenmsb  \scriptfont\msbfam=\sevenmsb
  \scriptscriptfont\msbfam=\fivemsb

\def\hexnumber@#1{\ifcase#1 0\or1\or2\or3\or4\or5\or6\or7\or8\or9\or
	A\or B\or C\or D\or E\or F\fi }

\def\Bbb{\ifmmode\let\next\Bbb@\else
 \def\next{\errmessage{Use \string\Bbb\space only in math mode}}\fi\next}
\def\Bbb@#1{{\Bbb@@{#1}}}
\def\Bbb@@#1{\fam\msbfam#1}
\def\Z{{\Bbb Z}}   \def\Q{{\Bbb Q}}  \def\Fb{{\Bbb F}\,} \def\R{{\Bbb R}}
\def\C{{\Bbb C}}
\font\smit=cmmi7

\rightline{October, 1995}
\bigskip\bigskip
\centerline{{\bf \huge  Kac-Peterson, Perron-Frobenius, and the}}
\bigskip\centerline{{\bf \huge Classification of Conformal Field Theories}}
\bigskip \bigskip   \centerline{Terry Gannon\footnote{$^{\dag}$}{{\little
Permanent address as of September 1996: Dept of Math, York University,
North York, Canada M3J 1P3}}}
\centerline{{\it Max-Planck-Institut f\"ur Mathematik, Bonn, D-53225}}
\bigskip\medskip
\centerline{Dedicated to the memory of Rudelle Hall, teacher and friend}
\bigskip\bigskip

{{\bf 1. Introduction: The classification of conformal field theories.}}
\quad Conformal field theories (CFTs) and related structures have
been of considerable value to mathematics, as for instance the work of
Witten has shown. This paper is concerned with their classification.
Fortunately, the problem has a simple expression in terms of the characters
of Kac-Moody algebras (see (1.2) below), and requires no prior knowledge of
CFT. Nevertheless, for reasons of  motivation, in the following paragraphs we
will sketch the definition of CFT.

Before discussing this background material, let us quickly state the actual
mathematical problem addressed in this paper. The characters of an affine
algebra at fixed level $k$ define in a natural way a unitary representation
of SL${}_2(\Z)$ (see equations (3.3) below). The ultimate classification
problem here is to find all matrices $M$ which commute with the matrices of
this representation, and which in addition obey relations (1.2b) and (1.2c)
-- such $M$ are called {\it physical invariants}. In this paper we address
the subproblem of finding all physical invariants which in addition
satisfy (1.3b), where $\s$ is the group of all symmetries of the (extended)
Coxeter-Dynkin diagram -- these $M$ we call {\it \ade s}.
Almost every physical invariant is expected to be a \ade.
In this paper we develop a program to find all of these for any affine
algebra, and apply it to explicitly find them for the algebra $A_r^{(1)}$.

The remainder of this introductory section is intended to
explain the motivation for this problem.
In the language of CFT (which will be touched on shortly), the
classification of these physical invariants is equivalent to the
classification of all possible {\it Wess-Zumino-Witten partition functions}.
There is, we shall see, a fairly natural cut of this classification
into two subproblems. One is to find all
possible {\it chiral algebras} (these are essentially vertex operator
algebras), and the other is to find all possible
automorphisms of the corresponding {\it fusion rings} (these encode the
tensor product structure of the algebra). In previous work
[11,12] we accomplished the second subproblem for the case where the chiral
algebra corresponds to an affine algebra; in this
paper we generalize those arguments to the case where the
chiral algebra is an extension of those by {\it simple currents} (see e.g.\
[6]). It is generally believed (for reasons given below) that `almost every'
chiral algebra finitely extending the chiral algebra of an affine algebra,
will be of this form, and
so this paper solves (for $A_r^{(1)}$) the second subproblem for what may be
termed its {\it generic} chiral extensions.

According to Segal [25] (see also the presentation in [15]), a
(two-dimensional) conformal field theory is a {\it representation}
of the category ${\cal C}$ whose objects are disjoint unions of parametrized
circles and whose morphisms are cobordisms --  i.e. it is a
functor ${\cal T}$ from ${\cal C}$ into the category of complex Hilbert
spaces and trace-class operators. There exists a Hilbert space $H$ such that
${\cal T}$ takes
$n$ circles to $H\otimes\cdots\otimes H$ ($n$ times). Sewing together
surfaces in ${\cal C}$ along boundary circles corresponds by ${\cal T}$ to
composing operators. The detailed definitions and axioms are not important
here, and would take us too far afield.

The data of a CFT decomposes into two {\it chiral halves}, related to the
fact that the conformal maps in $\C$ consist of analytic functions and
their complex conjugates. Of greatest interest are the rational conformal
field theories (RCFTs), defined by Segal using the notion of a {\it modular
functor}. The modular functor makes precise the constraints imposed on
each chiral half: the key property of an RCFT is that the chiral data
is labelled by a {\it finite} set (the primary fields of the theory).

${\cal T}$ will map the closed torus $\C/(\Z+\Z\,\tau)$ to a complex
number $\z(\tau)$; $\z$ is called the {\it partition function} for the
theory.  But different $\tau$ can correspond to the same torus; these $\tau$
are related by the modular group of the torus.
Thus the partition function $\z$ should be
modular invariant, i.e.\ invariant under the natural action of PSL${}_2(\Z)$
on the upper half complex plane.

Important examples of RCFTs are where the `chiral labels'
 are given by representations of a Kac-Moody algebra $X_r^{(1)}$ at some
fixed level $k\in\{1,2,3,\ldots\}$.
These are called Wess-Zumino-Witten (WZW) models\footnote{${}^1$}{{\little
WZW is
often used in the narrower sense of strings propagating on a group manifold,
so the term {\smit conformal}\ {\smit current}\ {\smit model}\ was proposed in
 [12] for the
more general case of interest here.}}. The partition function of a WZW model
will be of the form
$$\z(v)=\sum_{\mu,\nu}M_{\mu,\nu}\,\chi_\mu(v)\,\chi_\nu(v)^*\ ,\eqno(1.1)$$
where the parameter $v$ can be taken to lie in a Cartan subalgebra of
$X_r^{(1)}$. The sum in (1.1) is over all highest weights $P_+(X_r^{(1)},k)$;
one of these weights, denoted $k\Lambda_0$, is distinguished. This differs
{} the partition function $\z(\tau)$ discussed earlier, only by depending
on more variables.
There is a natural action $v\mapsto Av$ of SL$_2(\Z)$
on the Cartan subalgebra [17]. The function
in (1.1) obeys the following conditions:
$$\eqalignno{&\z(Av)=\,\z(v)\quad {\rm for\ all}\ A\in {\rm SL}_2(\Z)\
;&(1.2a)\cr
&M_{\mu,\nu}\in\,\{0,1,2,\ldots\}\ ;&(1.2b)\cr
&M_{k\Lambda_0,k\Lambda_0}=\,1\ .&(1.2c)\cr}$$
Any such $\z$ or $M$ is called a {\it physical invariant}.

WZW models have been extensively studied because they are simple enough to
analyze, but complicated enough that the answers should be interesting
and hopefully characteristic of more general RCFTs. They are generally
regarded as building blocks, via the Goddard-Kent-Olive coset construction,
for perhaps all other RCFTs.

One of the few remaining fundamental questions of WZW models is their
classification. Because of Segal's sewing axiom, the higher genus behaviour
of a CFT can be determined in principle from that of lower genus. In
particular, a RCFT is uniquely determined by its two chiral algebras (which
in most cases are taken to be isomorphic),  the operator product structure
coefficients (obtained in principle from ${\cal T}$ by selecting a disk
with two punctures), and the partition function.

In this paper we address the classification of possible WZW partition
functions (i.e.\ physical invariants). The first such result was for
$A_1^{(1)}$, for all
$k$ [4]. It was found that the set of solutions to (1.2) for $A_1^{(1)}$
fall into the mysterious A-D-E pattern (see also [26]). An explanation
has recently been announced by Ocneanu [21], using subfactor theory and
path algebras on graphs. All physical invariants are also known for
$A_2^{(1)}$ [10]. For it, no connection with A-D-E is known,
but several unexplained coincidences have appeared (see e.g.\ [22]) between
the $A_2^{(1)}$ classification and the Jacobians of Fermat curves.
Zuber [28] and collaborators have explored using generalized Coxeter graphs to
reinterpret and extend some of these observations.
Classifying RCFTs is interesting in its own right, but what makes it more
intriguing is the desire to understand and if possible generalize these
apparent patterns.

Unfortunately physical invariants have resisted
extensive attempts at their classification; only for $A_1^{(1)}$ [4],
$A_2^{(1)}$ [10], and $(A_1\oplus A_1)^{(1)}$ [9]\footnote{${}^2$}{
{\little However for this latter algebra an additional constraint beyond
(1.2), involving the
Knizhnik-Zamolodchikov equation, was assumed.}} has the classification
been attained at all levels $k$. However there has been recent
progress [11,12] toward the solution of this problem, and this paper
takes us one step closer to this goal.

Let $\s$ denote the group of all symmetries of the (extended)
Dynkin diagram of $X^{(1)}_r$. Any $A\in\s$ will induce
a permutation $\la\mapsto A\la$ of the level $k$ weights of $X_r^{(1)}$, by
the action of $A$ on the Dynkin labels. Write $\s\la$ for the orbit of $\la$
by $\s$. The $A\in\s$ which fix the extended node are called {\it
conjugations}; some of the remainder (defined in section 3 below) are called
{\it simple currents}. It is easy to verify (see (3.5)
below) that the modular behaviour of $\chi_{A\la}$ is closely related to
that of $\chi_\la$, for any symmetry $A\in\s$. So it is not surprising
that these can be used to obtain new physical invariants from
old ones [2]. Indeed it seems that most
physical invariants can be obtained in this way
{}from the identity matrix physical invariant $M=I$ -- such physical
invariants are called {\it simple current invariants} (and their
conjugations). See (1.4) below.

As can be seen from (1.2c), as well as (3.4c) below, the weight $k\Lambda_0$
has special significance. A reasonable division of this classification
problem into two subproblems is, on the one hand, to consider all possible
values $M_{k\Lambda_0,\mu}$, $M_{\mu,k\Lambda_0}$ -- these are severely
constrained [10] -- and on the other hand to find all physical invariants
$M$ which realize each of these possible choices for $M_{k\Lambda_0,\mu}$,
$M_{\mu,k\Lambda_0}$. This is a restatement of the two
subproblems mentioned in the third paragraph. In [11,12] we find all
possible physical invariants satisfying the additional constraint
$$M_{\mu,k\Lambda_0}\ne 0\ {\rm or}\ M_{k\Lambda_0,\mu}\ne 0\quad
\Longrightarrow \quad \mu=k\Lambda_0\ .\eqno(1.3a)$$
These are called {\it automorphism invariants}. In this paper we generalize
those arguments to find all physical invariants obeying instead the
more general constraint
$$M_{\mu,k\Lambda_0}\ne 0\ {\rm or}\ M_{k\Lambda_0,\mu}\ne 0\quad
\Longrightarrow \quad \mu\in \s(k\Lambda_0)\ .\eqno(1.3b)$$
We call these ${\cal ADE}_7${\it -type invariants}. For example, in the
$A_1^{(1)}$ classification [4], these consist of
the physical invariants called ${\cal A}_\ell$ and ${\cal D}_\ell$, along
with the exceptional ${\cal E}_7$. Based on the known classifications
(e.g.\ [4,10,9,5]), together with various computer checks in the
literature,
it is reasonable to conjecture that almost all physical invariants are
\ade s. For example, for all but a small number of
$A_r^{(1)}$, we expect all physical invariants for each $k\ne r-1,r+1,r+3$
to obey (1.3b). The \ade s are interesting also
because they include exceptional physical invariants (like ${\cal E}_7$ for
$A_1^{(1)}$) which are notoriously difficult to obtain by standard
constructions.

This paper is concerned with the classification of all \ade s.
We reduce the problem to the mechanical albeit tedious
task of computing q-dimensions and tensor product coefficients. We complete
the  classification for the case of greatest interest: $A_r^{(1)}$. Up to
conjugations, we find only 8 exceptional \ade s.
This is a significant step towards the classification of all
WZW partition functions for the unitary algebras. The final step in that
classification,
namely solving the various constraints for $M_{\mu,k\Lambda_0},M_{k\Lambda_0
,\mu}$, will not be  addressed here.

Some of the arguments in this paper are based on those for the automorphism
invariant classification [11,12], as well as older classifications
[9,10], but several
new complications arise here. The main  tools we use are the Kac-Peterson
formula (3.4d) -- which permits us to exploit the well-understood
representation  theory of finite-dimensional Lie algebras -- and the
Perron-Frobenius spectral theory for  non-negative matrices.

A somewhat related problem is [19] to classify all physical invariants
which for all weights $\mu,\nu$ obey the constraint
$$M_{\mu,\nu}\ne 0\Rightarrow \nu\in \s_{sc}\mu\ ,\eqno(1.4)$$
where $\s_{sc}$ is the subgroup of $\s$ consisting of all simple currents.
These are called {\it simple current invariants}; they are a special case
of the \ade s considered here. Their classification
has been accomplished for all RCFT, subject to a certain constraint on
the modular $S$ matrix (3.3c) [19] -- it is found that there are no
exceptional invariants of this form. Though this is clearly a major result,
(1.4) is sufficiently stronger than (1.3b) that the arguments in [19] are not
useful in our context.

In section 2 below we list all \ade s for $A_r^{(1)}$.
Section 3 establishes the basic results
we need, and section 4 specializes to $A_r^{(1)}$ and outlines the argument
for classifying all $A_r^{(1)}$ \ade s.
The problem reduces to some q-dimension calculations and computing some
tensor product coefficients, which we do in sections 5 and 6 respectively.
This completes the classification for almost all levels $k$ of $A_r^{(1)}$;
the finitely many trouble-making pairs $(r,k)$ are explicitly handled in
section 7.

\bigskip\bigskip{{\bf 2. The \ade s of
$A_r^{(1)}$.}}\quad
In this section  we explicitly list all of the \ade s
of $g=A_r^{(1)}$. The proof that this list
is complete  will be accomplished in the later sections. In the following
section we will motivate and generalize many of the definitions made here;
our purpose here is merely to state Theorem 2.1.

Fix the rank $r$ and level $k$, and define $\r=r+1$ and $\k=k+\r$.
The level $k$ highest weights of $A_r^{(1)}$ constitute the set
$P_{+}$ of
$\r$-tuples $\la=(\la_0,\ldots,\la_r)$ of non-negative integers $\la_i$
obeying $\sum_{i=0}^r\la_i=k$. The extended Coxeter-Dynkin diagram of
$A_r^{(1)}$ is a circle with $\r$ nodes, which
we label counterclockwise 0 to $r$ -- 0 is called the {\it extended} node.
Its $2\r$ symmetries (only $\r$, if $r=1$)
form the dihedral group $\s$; it
is generated by an order 2 symmetry fixing 0 (the {\it conjugation} $C$),
and an order $\r$ rotation taking $i$ to $i+1$ (the {\it simple current} $J$).
This group acts on $P_+$ by permuting the indices of the weight:
$$\eqalignno{C\la=&(\la_0,\la_r,\la_{r-1},\ldots,\la_1)\ ,&(2.1a)\cr
J\la=&(\la_r,\la_0,\la_1,\ldots,\la_{r-1})\ .&(2.1b)\cr}$$
A convenient quantity we will often use is the $\r$-ality $t$ defined by
$$t(\la)\eqde \sum_{j=1}^r j\la_j\ .\eqno(2.2)$$

Together with Lemmas 3.1, 3.2 and 3.3, the following theorem is the main
result of this paper. The \ade s named in Theorem 2.1
are defined in equations (2.3a), (2.4), (2.5), and (2.7) below.

\medskip{{\smcap Theorem 2.1.}} \quad{\it The complete list of \ade s
for $A_{r}^{(1)}$ at level $k$ is:}

\item{$\bullet$} {\it for all $r,k\ge 1$, $d$ dividing $\r$ and satisfying
(2.3b), and} $c=0,1$:\quad $C^c\cdot I[\J_d]$\ ;

\item{$\bullet$} {\it for} $(r,k)\in\{(1,16),\,(3,8),\,(4,5),\,(7,4)\}$:
\quad $\E^{(r,k)}$;

\item{$\bullet$} {\it for $(r,k)\in\{(2,9),\,(8,3)\}$}:
\quad\qquad\qquad $\E^{(r,k)}$ {\it and} $C\cdot \E^{(r,k)}$;

\item{$\bullet$} {\it for $(r,k)=(15,2)$}:
\quad\qquad\qquad\qquad $\E^{(15,2)}$, ${1\over 2}\,I[\J_4]\cdot\E^{(15,2)}$
{\it and} $C\cdot \E^{(15,2)}$.

\medskip Next, we explicitly define these \ade s.

Denote by $\J_d$ the subgroup of $\s$ generated by $J^d$, when $d$ divides
$\r$. Each such subgroup can be used to construct a \ade.
In particular, put ${k}'=\k$ if both $k$ and $\r$ are odd, otherwise
put ${k}'=k$. Define [23]
$$I[\J_d]_{\la,\mu}=\sum_{j=1}^{\r/d}\delta^{\r/d}(t(\la)
+d{}j{k}'/2)\ \delta_{\mu,J^{d{}j}\la}\ ,\eqno(2.3a)$$
where $\delta^y(x)=1$ or 0 depending, respectively, on whether or not
$x/y\in\Z$. Then $I[\J_d]$ will be a physical invariant iff [23]
$$ {k}'d\equiv 0\qquad ({\rm mod}\ 2)\ .\eqno(2.3b)$$
This can be readily proven using (4.1b) and (4.1c) below.

These $I[\J_d]$ were first explicitly given in [7], though some appeared
earlier in [2]. Equation (2.3a) extends naturally to any $X_r^{(1)}$ (see
[23]). Note that $d=\r$ always satisfies (2.3b); it gives $I[\J_{\r}]=I$, the
identity matrix. Incidently, for each divisor $d$
of $\r$, there is a Lie group $G_d$ whose simply-connected covering group
is $G_1\eqde{\rm SU}_{\r}$, and which obeys $\|G_1/G_d\|=d$.
For example, for $\r=2$ $G_{2}={\rm SO}_{3}$. The existence of $I[\J_d]$
is intimately connected to that of $G_{\r/d}$ [7].

The conjugation $C$ defines another \ade\ (see
(4.1d), (4.1e)), which we will also denote by $C$:
$$C_{\la,\mu}=\delta_{\mu,C\la}\ ,\eqno(2.4)$$
where $\delta$ denotes the Kronecker delta. Moreover, the matrix product $C
\cdot M$ of $C$ with any other \ade\ $M$ will also be a
\ade, and $C^2=I$.

In addition, there are a number of other \ade s,
called ${\cal E}_7$-{\it type exceptionals}. It is slightly more convenient
to express these in terms of characters rather than their coefficient
matrices $M$.
It suffices to give the relevant subgroup $\J_d$, as well as  the characters
 with the exceptional behaviour (the remaining characters combine exactly
 as in $I[\J_d]$). To help explain our
notation, we will write out in full the two simplest such exceptionals:
$$\eqalignno{\E^{(1,16)}=&\,|\chi_{16,0}+\chi_{0,16}|^2+|\chi_{12,4}+
\chi_{4,12}|^2+|\chi_{10,6}+\chi_{6,10}|^2&\cr
&\,+(\chi_{14,2}+\chi_{2,14})\,\chi_{8,8}^*+\chi_{8,8}\,(\chi_{14,2}+
\chi_{2,14} )^*+|\chi_{8,8}|^2\ ;&(2.5a)
\cr\E^{(2,9)}=&\,|\chi_{900}+\chi_{090}+\chi_{009}|^2+|\chi_{522}+\chi_{252}+
\chi_{225}|^2+|\chi_{603}+\chi_{360}+\chi_{036}|^2&\cr
&+|\chi_{630}+\chi_{063}+\chi_{306}|^2+|\chi_{144}+\chi_{414}+\chi_{441}|^2+
2|\chi_{333}|^2&\cr
&+(\chi_{711}+\chi_{171}+\chi_{117})\,\chi_{333}^*+\chi_{333}\,(\chi_{711}+
\chi_{171}+\chi_{117})^*\ .&(2.5b)\cr}$$
Here and elsewhere, we label a weight by its Dynkin labels.

For convenience write
$$\langle\chi_\la\rangle_d\eqde \sum_{\mu\in\J_d\la}\chi_\mu\ .\eqno(2.6a)$$
Also, write ``$a* b$'' as short-hand for ``$ab^*+ba^*$''.
Note that we may capture all the information in (2.5a) by stating
$d=1$, and giving the `exceptional' terms:
$$\langle\chi_{14,2}\rangle_1*\chi_{8,8}+|\chi_{8,8}|^2\ . \eqno(2.6b)$$
The remaining terms in (2.5a) are exactly as in $I[\J_1]$.
Similarly, (2.5b) can be summarized by stating $d=1$, and
giving the exceptional terms
$$2|\chi_{333}|^2+\langle\chi_{711}\rangle_1*\chi_{333}\ .\eqno(2.6c)$$

The remaining ${\cal E}_7$-type exceptionals are expressed in this way as:
$$\eqalignno{\E^{(3,8)}:\quad d=1;\quad& 2|\chi_{2222}|^2
+(\langle\chi_{5012}\rangle_1+\langle\chi_{5210}\rangle_1)*
\chi_{2222}&\cr&+\langle\chi_{6101}\rangle_1*\langle\chi_{4040}
\rangle_1+|\langle\chi_{4040}\rangle_1|^2\ ;&(2.7a)\cr
\E^{(4,5)}:\quad d=1;\quad& \langle\chi_{31001}\rangle_1*\chi_{11111}
+4|\chi_{11111}|^2\ ;&(2.7b)\cr
\E^{(7,4)}:\quad d=2;\quad& |\langle \chi_{20002000}\rangle_2|^2+|\langle
\chi_{02000200}\rangle_2|^2+\langle\chi_{21000001}\rangle_2*\langle
\chi_{20002000}\rangle_2&(2.7c)\cr
&+\langle\chi_{12100000}\rangle_2*\langle\chi_{02000200}\rangle_2 +(\langle
\chi_{12000010}\rangle_2+\langle\chi_{10100002}\rangle_2)*\chi_{01010101}&\cr
&(\langle\chi_{21010000}\rangle_2+\langle\chi_{20000101}\rangle_2)*
\chi_{10101010}+2|\chi_{01010101}|^2+2|\chi_{10101010}|^2\bigr\};&\cr
\E^{(8,3)}:\quad d=3;\quad& \sum_{j=0}^2\bigl(\, 2|\chi_{J^j(\L_0+\L_3+\L_6)}
|^2+\langle\chi_{J^j(\L_2+\L_3+\L_4)}\rangle_3*\chi_{J^j(\L_0+\L_3+\L_6)}
\bigr)\ ;&(2.7d)\cr
\E^{(15,2)}:\quad d=8;\quad &\sum_{j=0}^7\bigl(\,|\chi_{J^j
(\Lambda_0+\Lambda_8)}|^2 +\langle\chi_{J^j(\Lambda_3+\Lambda_5)}
\rangle_8*\chi_{J^j(\Lambda_0+\Lambda_8)}\bigr)\ .&(2.7e)\cr}$$

$\E^{(1,16)}$ was first given in [4]; $\E^{(2,9)}$ in [20];
$\E^{(4,5)}$ in [24]; $\E^{(3,8)}$, $\E^{(7,4)}$, $\E^{(8,3)}$ in [8].
$\E^{(15,2)}$ is new but [8] obtained its projection: the matrix product
${1\over 2}\,I[\J_4]\cdot\E^{(15,2)}$, which has $d=4$ and the exceptional
terms
$$\sum_{j=0}^3\bigl(|\chi_{J^{2j}
(\Lambda_0+\Lambda_8)}|^2 +\langle\chi_{J^{2j}(\Lambda_3+\Lambda_5)}
\rangle_4*\chi_{J^{2j}(\Lambda_0+\Lambda_8)}\bigr)\ .\eqno(2.7f)$$

Note the symmetry $(r,k)\leftrightarrow(k-1,r+1)$ in the list of ranks and
levels of these exceptionals. This is not surprising, considering
the {\it rank-level duality} (see (4.2), (4.3) below) exhibited by the
Kac-Peterson $S$ and $T$ matrices.

\medskip{{\it Remarks 2.1.}} Note that the matrices $M$ of all \ade s
here are symmetric: $M=M^T$. This is not always true for other
$g$ [13]. There are some redundancies in the list in Theorem 2.1. When
$r=1$ or $k\le 2$, take $c=0$ only -- this is because there $C$ will equal
one of the $I[\J_d]$. Likewise, $C\cdot I[\J_d]=I[\J_d]$
for $(r,k,d)\in\{(2,3,1),\,(2,6,1),\,(4,5,1),\,(5,3,2)\}$. The final
redundancy is $I[\J_1]=I[\J_2]$ for $\r=k=2$.

\bigskip\bigskip
{{\bf 3. Cyclotomy, Kac-Peterson and Perron-Frobenius.}}\quad
 In this section we establish the fundamental lemmas which define our
 program to classify all \ade s. We will state
and prove them for any $g=X_r^{(1)}$ -- indeed they continue to hold for
any RCFT. The notation used here is standard; see e.g.\ [16] for
more details. We will quickly review the basic facts, before heading
into the statement and proof of Lemmas 3.1, 3.2, 3.3.

Let $g$ be the non-twisted affine algebra $X_r^{(1)}$ derived from the
finite-dimensional algebra $\bar{g}=X_r$.
Let $L(\la)$ denote any irreducible integrable highest weight $g$-module,
and let $\chi_\la$ be its normalized character
with respect to a Cartan subalgebra $h=\bar{h}\oplus{\Bbb C}c\oplus {\Bbb C}
d$ ($\bar{h}$ is a Cartan subalgebra of $\bar{g}$, and $c$ the canonical
central element, $d$ a derivation, of $g$).

Let $\L_0,\ldots,\L_r\in h^*$ denote the fundamental weights of $g$. Then
the highest weight $\la\in P_+(g,k)$ of $L(\la)$ can be taken to lie in
$$P_+\,\eqde \bigl\{\sum_{j=0}^r\la_j\L_j\,
|\,\la_j\in \Z,\ \la_j\ge 0,\ \sum_{j=0}^ra^\vee_j\la_j=k\bigr\}\ ,
\eqno(3.1)$$
where $k$ is a positive integer called the level, and the positive integers
$a_j^\vee$ are the co-labels of $g$. Write $h^\vee=\sum_i a^\vee_i$, and
$\bar{k}=k+h^\vee$. The Weyl vector
is $\rho=\sum_{i=0}^r \L_i$. For convenience we will write $\L^0$ for $k\L_0$.
Note that the projection
$$\la\mapsto \overline{\la}\,\eqde \sum_{i=1}^r\la_i\bar{\L}_i,\qquad
\bar{\L}_i\eqde \L_i-a_i^\vee\L_0$$
produces a highest
weight for the underlying algebra $\overline{g}=X_r$. This projection
is orthogonal with respect to the invariant bilinear form $(-|-)$ (which
we take to be normalized so that long roots have norm 2); in fact
$(\la|\la)=(\bar{\la}|\bar{\la})$.

Let $Q^\vee$ denote the coroot lattice, and $\overline{W}$ the Weyl
group of $\bar{g}$.
The affine Weyl group $W$ is isomorphic to the semi-direct product $T\cdot
\overline{W}$, where $T$ consists of the translations $t_\al$, $\al\in
Q^\vee$, defined on $h^*$ (mod $\C\,\delta)$ by
$$t_\al\L_i\equiv \L_i+a_i^\vee\al\qquad({\rm mod}\ \C\,\delta)
\eqno(3.2)$$
($\delta$ here is an imaginary root of $g$). This is a central observation in
the representation
theory of affine Kac-Moody algebras. It permits an expression -- the Weyl-Kac
character formula -- for the character
$\chi_\la$ in terms of theta functions. This implies [17] that $\chi_\la$
will be a modular function: in particular, we may regard $\chi_\la$ as a
function from $h$ to $\C$; coordinatizing $h$ in the usual
way (i.e.\ $2\pi i(z-\tau\,d+u\,c)\in h$, where $z\in \bar{h}$,
$\tau,u\in\C$), we obtain
$$\eqalignno{\chi_\la(\tau+1,z,u)&=\sum_{\mu\in P_+}T_{\la,\mu}
\,\chi_\mu(\tau,z,u)\ ,&(3.3a)\cr
T_{\la,\mu}&\eqde \exp[\pi i\{{(\la+\rho|\la+\rho)\over \overline{k}}-
{(\rho|\rho)\over h^\vee}\}]\,\delta_{\la,\mu}\ ;&(3.3b)\cr
\chi_\la({-1\over\tau},{z\over \tau},u-{(z|z)\over 2\tau})&=\sum_{\mu\in
P_+} S_{\la,\mu}\,\chi_\mu(\tau,z,u)\ ,&(3.3c)\cr
S_{\la,\mu}&\eqde s
\sum_{w\in \overline{W}}{\rm det}(w)\,\exp[-2\pi i{(w(\la+\rho)|\mu+\rho)\over
\overline{k}}]\ ;&(3.3d)\cr}$$
where in (3.3d) the normalization $s$ is
$$s=i^{\|\overline{\Delta}_+\|}\,\bar{k}^{-r/2}\,\|Q^{\vee\ *}/
Q^\vee\|^{-{1\over 2}}\ .$$
Here, $\|\overline{\Delta}_+\|$ denotes the number of positive roots of
$\overline{g}$, and the weight lattice $Q^{\vee\ *}$ is the dual lattice of
$Q^\vee$. Together, (3.3a),(3.3c) define the transformation properties
of $\chi_\la$ with respect to SL$_2(\Z)$.

These {\it Kac-Peterson} matrices $S$ and $T$ have some special properties.
They are unitary and symmetric. From the Weyl denominator formula we get,
for any $\ell=0,\ldots,k/(h^\vee-1)$,
$$S_{\ell\bar{\rho},\la}=|s|\,\prod_{\bar{\al}>0} 2\sin\bigl( \pi
{(\bar{\la}+\bar{\rho}|\bar{\al})\over \bar{k}/(\ell+1)}\bigr)\ ,\eqno(3.4a)$$
where by $\ell\bar{\rho}$ in (3.4a) we mean the weight
$\ell\bar{\rho}+(k-\ell(h^\vee-1))
\L_0$, and where the product is over the positive roots $\bar{\al}\in
\overline{\Delta}_+$ of $\bar{g}$.
Usually we will take $\ell=0$ in (3.4a).
This implies the following expression for the {\it q-dimensions:}
$$\D(\la)\eqde{S_{\la,\L^0}\over S_{\L^0,\L^0}}=\prod_{\bar{\al}>0}
{\sin(\pi \,(\bar{\la}
+{\bar\rho} |\bar{\al})/\k)\over \sin(\pi\,(\bar{\rho}|\bar{\al})/\k)}\ .
\eqno(3.4b)$$
{}From (3.4a) one can show that
$$S_{\la,\L^0}\ge S_{\L^0,\L^0}>0,\quad\forall \la\in P_+\
.\eqno(3.4c)$$
$S$ also satisfies the important equation [17]
$${S_{\la,\mu}\over S_{\L^0,\mu}}=\overline{ch}_{\overline{\la}}(-2\pi i
{\overline{\mu}+\overline{\rho}\over\k})\ ,\eqno(3.4d)$$
where $\overline{ch}_{\overline{\la}}$ is the Weyl character of the
$\overline{g}$-module $\overline{L}(\bar{\la})$. Equation (3.4d) has many
consequences, one of which has to do with the {\it fusion
coefficients} of $g$. These can be taken to be defined by Verlinde's
formula:
$$N_{\la,\mu}^\nu\eqde \sum_{\ga\in P_+}S_{\la,\ga}\,{S_{\mu,\ga}\over
S_{\L^0,\ga}}\,S_{\nu,\ga}^*\ .\eqno(3.4e)$$
Fusion coefficients have an algebraic interpretation in terms of the
tensor products of representations of e.g.\ Hecke  algebras and quantum groups
at roots of unity, as well as a geometric
interpretation involving moduli spaces of principle bundles over projective
curves. In the language of RCFTs, they give the dimensions of the spaces
of conformal blocks. The only relevant point here is that, because
of (3.4d), they can be computed in terms of the tensor product
multiplicities mult$_{\overline{\la}\otimes\overline{\mu}}(\overline{\nu})$
in $\bar{g}$ [27,16], and hence its weight multiplicities
$m_{\overline{\la}}(\overline{\mu})\eqde {\rm dim}\,\overline{L}(\overline{\la}
)_{\overline{\mu}}$:
$$N_{\la,\mu}^\nu=\sum_{w\in{W}}{\rm det}(w)\,m_{\overline{\mu}}
(\overline{w.\nu}-\overline{\la})\ ,\eqno(3.4f)$$
where $w.\ga\eqde w(\ga+\rho)- \rho$ (compare the Racah-Speiser algorithm
for computing mult$_{\overline{\la}\otimes\overline{\mu}}(\overline{\nu})$).

The symmetries of the (extended) Coxeter-Dynkin diagram of $g$ define the
group $\s$. These play a major role in this paper. Those fixing
the extended node are called {\it conjugations}.  Another subgroup is
$\s_{sc}=W_0^+$, defined as follows. Let $T_0$
denote the set of all translations $t_\al$ in (3.2) with $\al\in P^\vee$,
where $P^\vee$ is the co-weight lattice. Define [18]
$$\s_{sc}=\{J\in T_0\cdot \overline{W}\,|\,J(\Delta_+)=\Delta_+\,\}\ ,$$
where $\Delta_+$ are the positive roots of $g$. $\s_{sc}$ stabilizes
$\Pi^\vee$, and defines a normal
subgroup of $\s$ isomorphic to $Q^\vee/P^\vee$. Its elements are called
{\it simple currents}. Both conjugations and simple currents act on
$P_+$  by permuting the Dynkin labels, and together they generate $\s$.

All conjugations commute with $S$ and $T$ and fix $\L^0$. The most important
conjugation is called {\it charge conjugation}: it takes each weight $\la$ to
$C\la={}^t\la$, the weight contragredient to $\la$. It obeys the important
relation
$$C=S^2\ .\eqno(3.5a)$$

The primary reason for the importance of $\s_{sc}$ is: let $J=t_\al w\in
\s_{sc}$, then [18]
$$S_{J^a\la,J^b\mu}=\exp[2\pi i\,(a\,Q_J(\mu)+b\,Q_J(\la)+ab\,Q_J(J\L^0))]\,
S_{\la,\mu}\ ,\eqno(3.5b)$$
where by $J^a$ we mean the $a$-fold composition $J\circ\cdots\circ J$, and
where $Q_J(\mu)=-(\bar{\mu}|\al)$. The matrix $T$ also behaves similarly
under $\s_{sc}$:
$$\exp[2\pi i \,Q_J(\la)]=T_{\la,\la}T^*_{J\la,J\la}T_{J\L^0,J\L^0}T^*_{\L^0,
\L^0}\ .\eqno(3.5c)$$

\medskip{{\it Definition 3.1.}}\quad By a {\it positive invariant}
for a given algebra $X_r^{(1)}$ and level $k$ we mean a matrix $M$ commuting
with the
corresponding Kac-Peterson matrices $S$ and $T$, with the additional
property that each $M_{\la,\mu}\ge 0$. By a {\it physical invariant} we
mean a positive invariant with each $M_{\la,\mu}\in\Z$, and obeying (1.2c).
By a {\it \ade}\ we mean a
physical invariant $M$ satisfying (1.3b).\medskip

For example, any conjugation defines a \ade. Simple
currents can also be used to construct them (see e.g.\ [2,23]) -- an
example is (2.3a). Any physical invariant not constructable in these
standard ways out of simple currents and conjugations  is called an
{\it exceptional invariant}, and if it is in addition a \ade,
we shall call it an ${\cal E}_7$-type exceptional (by analogy
with the A-D-E classification in [4]).

The condition $TM=MT$ is equivalent to the `selection rule'
$$M_{\la,\mu}\ne 0\Rightarrow (\la+\rho\,|\,\la+\rho)\equiv(\mu+\rho\,|\,
\mu+\rho)\ \quad({\rm mod}\ 2\bar{k})\ .\eqno(3.6a)$$
The other commutation condition, namely
$$S\,M=M\,S\ ,\eqno(3.6b)$$
or equivalently  (since $S$ is unitary)
$$S\,M\,S^{\dag}=M\ ,\eqno(3.6c)$$
is much more subtle and interesting, and we will begin to explore its
consequences in this section. Equations (3.6a) and (3.6b) are equivalent
to the modular invariance condition (1.2a).

For a positive invariant $M$, define
$$\eqalignno{\J_L(M)=&\,\{J\in \s_{sc}\,|\,M_{J\L^0,\L^0}\ne 0\}\ ;&(3.7a)\cr
\p_L(M)=&\,\{\la\in P_+\,|\,\exists \mu\in P_+\ {\rm such\ that}\
M_{\la,\mu} \ne 0\}\ ;&(3.7b)\cr}$$
and define $\J_R(M)$ and $\p_R(M)$ similarly (using the other subscript of
$M$). Call $\la\in P_+$ a {\it fixed point} of $\J\subset \s_{sc}$ if
$\|\J\la\|<\|\J\|$. Let $\F(\J)$ denote the set of all fixed points of $\J$.
For any $\J\subset \s_{sc}$, define
$$\p(\J)\eqde\{\la\in P_+\,|\,Q_J(\la)\equiv 0\ ({\rm mod}\ 1)\ \quad
\forall J\in\J\}\ .\eqno(3.7c)$$

The remainder of this section is devoted to the statement and proof of
the basic lemmas we will need. Exactly how to use these will be addressed
in the following section.

Our first lemma is an easy consequence of (3.4c) and (3.5b).
It tells us how $\J_L(M)$ and $\J_R(M)$
influence all other values of $M$.

\medskip{{\smcap Lemma 3.1}}.\quad(a) {\it Let $M$ be any physical invariant,
and $J,J'\in \s_{sc}$. Then the following statements are equivalent:}

\qquad\qquad\qquad{(i)} $M_{J\L^0,J'\L^0}\ne 0$;

\qquad\qquad\qquad{(ii)} $M_{J\L^0,J'\L^0}=1$;

\qquad\qquad\qquad{(iii)} {\it for any $\la,\mu\in P_+$, if
$M_{\la,\mu}\ne 0$ then $Q_J(\la)\equiv Q_{J'}(\mu)$ (mod 1)};

\qquad\qquad\qquad{(iv)} $M_{J\la,J'\mu}=M_{\la,\mu}$ {\it for all} $\la,\mu
\in P_+$.

\item{{(b)}} \quad {\it Let $M$ be any positive invariant satisfying
$$M_{\L^0,\mu}=\sum_{J\in\J_R}\delta_{\mu,J\L^0},\qquad M_{\la,\L^0}=
\sum_{J\in\J_L}\delta_{\la,J\L^0}\ ,\eqno(3.7d)$$
for some $\J_L,\J_R$. Then}

\qquad\qquad\qquad(i)  $M_{J\la,J'\mu}=M_{\la,\mu}$ {\it for all $\la,\mu
\in P_+$ and all} $J\in\J_L$, $J'\in\J_R$.

\qquad\qquad\qquad(ii) $\J_L$ {\it and $\J_R$ are groups and} $\|\J_L\|
=\|\J_R\|$.

\qquad\qquad\qquad(iii) $\p_L(M)=\p(\J_L)$ {\it and} $\p_R(M)=\p(\J_R)$.

\medskip{{\it Proof.}} \quad(a)  Note that
$$M_{J\la,J'\mu}=\sum_{\be,\ga}S_{J\la,\be}\,M_{\be,\ga}\,S^*_{\ga,J'\mu}
=\sum_{\be,\ga}\exp[2\pi i(Q_J(\be)-Q_{J'}(\ga))]\,S_{\la,\be}\,M_{\be,\ga}
\,S^*_{\ga,\mu}\ .\eqno(3.8a)$$
Applying this to $\la=\mu=\L^0$, and using (3.4c), we get that
$|M_{J\L^0,J'\L^0}|\le| M_{\L^0,\L^0}|$ with equality iff
the condition (iii) holds.  Thus for any physical invariant $M$,
(i) $\Rightarrow$ (ii) $\Rightarrow$ (iii). Statement (iii) implies (iv)
by (3.8a), and (iv) $\Rightarrow$ (i) is trivial.

\smallskip{(b)}  The argument from (3.8a) with $\la=\mu=\L^0$ and
$J'=id.$ tells us that $J\in\J_L$ if $Q_J(\beta)\in\Z$ for all $\beta\in
\p_L(M)$ -- this implies, again from (3.8a), that
$M_{J\la,\mu}=M_{\la,\mu}$
for all $J\in\J_L$, for all $\la,\mu\in P_+$. Similarly for $\J_R$. Hence
$\J_L,\J_R$ are groups.
The relation $\|\J_L\|=\|\J_R\|$  comes from the calculation
$$S_{\L^0,\L^0}\,\|\J_L\|=\sum_\ga S_{\L^0,\ga}\,M_{\ga,\L^0}=\sum_\ga
M_{\L^0,\ga}\,S_{\ga,\L^0}=S_{\L^0,\L^0}\,\|\J_R\|\ .\eqno(3.8b)$$
In addition, $(SM)_{\la,\L^0}=(MS)_{\la,\L^0}$ gives us
$$S_{\la,\L^0}\sum_{J\in\J_L}\exp[2\pi i \,Q_J(\la)]=\sum_\ga M_{\la,\ga}\,
S_{\ga,\L^0}\ .\eqno(3.8c)$$
The r.h.s.\ is $>0$ iff $\la\in \p_L(M)$. The l.h.s.\ is $>0$ iff $Q_{J}(\la)
\in\Z$ for all $J\in \J_L$, since $\J_L$ is a group. \quad
\QED\medskip

Of course by Lemma 3.1(a), any \ade\ $M$ will obey
(3.7d) with $\J_L=\J_L(M)$ and $\J_R=\J_R(M)$. The converse however is
false. We will state and prove the remaining results in this section for
positive invariants obeying (3.7d), even though our primary interest is in
\ade s.\medskip

{{\it Definition 3.2.}} \quad For a given positive invariant $M$ obeying
(3.7d), call the pair $(\la,\mu)\in\p_L
(M)\times\p_R(M)$ $M$-{\it monogomous} if for all $\nu\in P_+$,  both
$$M_{\la,\nu}\ne 0 \Rightarrow \nu\in\J_R\mu\ ,\qquad{\rm and}\qquad
M_{\nu,\mu}\ne 0 \Rightarrow \nu\in\J_L\la\ .$$
In this case we also say $\la$ (resp.\ $\mu$) is {\it right}-(resp.\
{\it left}-)$M$-{\it monogomous}.\medskip

Note that (3.7d) says that $(\L^0,\L^0)$ is $M$-monogomous. For another
example,  every $\la\in\p_L(M)$ for $M=I(\J_d)$ (see (2.3a)) is
right-$M$-monogomous. We will find that $M$-monogomous pairs are the basic
building blocks of \ade s. When $(\la,\mu)$ is
$M$-monogomous, the value of $M_{\la,\mu}$ is given by (3.9) below. Also,
Lemma 3.2(b) below implies that whenever $\la$ and $\mu$ are not fixed points
of $\J_L$ and $\J_R$ respectively, then $M_{\la,\mu}\ne 0$, for some
\ade\ $M$, means $(\la,\mu)$ must be $M$-monogomous.

Another tool we will use is the Perron-Frobenius theory of non-negative
matrices --
see e.g.\ [14]. By a non-negative matrix $B$ is meant one whose entries
are all non-negative real numbers. Such a matrix has an eigenvalue $r(B)\ge 0$
with the property that $r(B)$ is at least as large as the modulus of any
other eigenvalue of $B$. An eigenvector corresponding to $r(B)$ is also
non-negative. For example, if $B$ is the $n\times n$ matrix satisfying
$B_{ij}=m$ for all $i,j$, then $r(B)=mn$. There are other properties of
non-negative matrices which
we will need below; we will state them as we use them. The next lemma
uses Perron-Frobenius to severely constrain the form $M$ can take.

For Lemma 3.2 and elsewhere, it is convenient to introduce the direct sum
decomposition
$$M=\oplus_j M_j=\left(\matrix{M_1&&0\cr &\ddots&\cr 0&&M_m\cr}\right)\ ,$$
where each $M_j$ is indecomposable (i.e.\ cannot be written as $M_j'\oplus
M_j''$). Let $\I(M_j)$ be the index set of
$M_j$. We will always take $M_1$ to be the unique one with $\L^0\in \I(M_1)$.

By an {\it irreducible matrix} in Lemma 3.2(a) below we mean a matrix which
cannot, under any simultaneous permutation of row and column indices, be
written in the form
$$\left(\matrix{A&B\cr 0&D\cr}\right)$$
for submatrices $A$, $B$,  $D$. Irreducible non-negative matrices have
special properties [14], as we shall see in the proof of Lemma 3.3 given
later.

\medskip{{\smcap Lemma 3.2}}.\quad (a) {\it Let $M$ be any {\it positive}
invariant satisfying (3.7d) for $\J=\J_L=\J_R$.
Then for each $i$, either $M_i=(0)$ or $r(M_i)=
\sqrt{r(M_iM_i^T)}=\|\J\|$. Moreover, each nonzero  $M_i$ is irreducible.}

(b) {\it  Let $M$ be any positive invariant satisfying (3.7d).
If $(\la,\mu)$ are $M$-monogomous, then
$$M_{J\la,J'\mu}={\|\J_L\|\over \sqrt{\|\J_L\la\|\, \|\J_R\,\mu\|}}\ ,\qquad
\forall J\in\J_L\ ,\ J'\in\J_R\ .\eqno(3.9)$$
Suppose $\be,\ga\in P_+$ satisfy $M_{\be,\ga}\ge\|\J_L\|/\sqrt{\|\J_L
\be\|\,\|\J_R\ga\|}$. Then $(\be,\ga)$ is $M$-monogomous.}

\medskip{\it Proof.}\quad (a) Consider $M(\ell)\eqde (M/\|\J\|
)^\ell$. The matrix $\|\J\|\,M(\ell)$ will also satisfy (3.7d). $M(\ell)$
will be a direct sum of $(M_i/\|\J\|)^\ell$. As
$\ell\rightarrow\infty$, $(M_i/\|\J\|)^\ell$ will tend identically to 0
if $r(M_i)<\|\J\|$, and will be unbounded if $r(M_i)>\|
\J\|$ -- both these follow for example from the Jordan canonical form
of $M_i$. That each $r(M_i)\le \|\J\|$ then follows from the (crude) bound
$${\rm max}_{\la,\mu} M'_{\la,\mu}\le \sum_{\la,\mu}M'_{\la,\mu}\le
{1\over (S_{\L^0,\L^0})^2}\sum_{\la,\mu}S_{\L^0,\la}\,M'_{\la,\mu}\,S_{\mu,
\L^0}={M'_{\L^0,\L^0}\over (S_{\L^0,\L^0})^2}\ ,\eqno(3.10)$$
valid for any positive invariant $M'$ (see (3.4c)).

The bound (3.10) also means that the sequence $\{M(\ell)\}_{\ell=1}^\infty$
will have a limit point $\widetilde{M}$, by Bolzano-Weierstrass.
$\|\J\|\,\widetilde{M}$ will be a positive invariant, and
will also satisfy (3.7d) for $\J=\J_L=\J_R$. By Lemma 3.1(b), this
means $\p_L(\widetilde{M})=\p_L(M)$, which forces $r(M_i)=\|\J\|$ whenever
$M_i\ne 0$.  That $r(M_iM_i^T)=\|\J\|^2$ follows by applying this result
to $M\,M^T/\|\J\|$.

Finally, to see that each $M_i\ne 0$ is irreducible, it suffices to show [14]
that both $M_i$ and $M_i^T$ have a strictly positive eigenvector corresponding
to eigenvalue $\|\J\|$. Let $v$
denote the vector with component $v_\mu=S_{\L^0,\mu}$ for each $\mu\in
\I(M_i)$. Then for each $\la\in \I(M_i)$,
$$\sum_{\mu\in\I(M_i)}(M_i)_{\la,\mu} \,v_\mu=\sum_{\mu\in P_+}
M_{\la,\mu}\, S_{\mu,\L^0}=\sum_{\mu\in P_+}S_{\la,\mu}\,M_{\mu,\L^0}=
\|\J\|\,S_{\la,\L^0}=\|\J\|\,v_\la\ ,$$
so $v$ is a positive eigenvector for $M_i$. The identical calculation and
conclusion holds for $M^T_i$.

\smallskip(b)\quad Put $m=\|\J_L\|=\|\J_R\|$. Let $\oplus_i
B_i$ be the direct sum decomposition of $MM^T/m$, where each $B_i$ is
indecomposable and $\L^0\in\I(B_1)$. $MM^T/m$ satisfies (3.7d) with $\J=\J_L$,
so $r(B_i)=m$ for all nonzero $B_i$.

First let us prove (3.9). Suppose $\la\in\I(B_i)$. Note that by Lemma 3.1(b),
$$(B_i)_{\la,\la}={\|\J_R\,\mu\|\over m}\,M^2_{\la,\mu}\ ,$$
so by Lemma 3.1(b) and the hypothesis that $(\la,\mu)$ is $M$-monogomous, we
get
$$m=r(B_i)={\|\J_R\,\mu\|\over m}\,M_{\la,\mu}^2\,\|\J_L\la\|\ .$$

Next, suppose $\be\in\I(B_j)$. Note that, by Lemma 3.1(b), $B_j\ge B_j'$
(component-wise), where
$$(B'_j)_{\nu\nu'}=\left\{\matrix{m/\|\J_L\,\beta\|&{\rm if}\ \nu,\nu'\in
\J_L\be\cr
0&{\rm otherwise}\cr}\right. \ .$$
Because $B_j$ is irreducible, we have [14] $r(B_j)\ge
r(B_j')$, with equality iff $B_j=B_j'$. However,
$r(B_j')=m$. Hence $B_j=B_j'$. A similar argument applied to $M^TM/m$ now
concludes the proof that $(\be,\ga)$ is $M$-monogomous.\quad\QED\medskip

{{\it Definition 3.3.}}\quad Given $\G\subseteq\p(\J)$, define
$$\G^{(1)}=\bigl\{\la\in\p(\J)\,|\,\forall \mu\in\p(\J),\ \mu\not\in\J\la,\
\exists\ga\in\G\ {\rm
satisfying}\ {S_{\ga,\la}\over S_{\L^0,\la}}\ne{S_{\ga,\mu}\over S_{\L^0,\mu}}
\,\bigr\}\ ,\eqno(3.11)$$
and $\G^{(n)}=(\G^{(n-1)})^{(1)}$. $\Gamma$ is called a {\it fusion-generator}
for $\J$ if $\G^{(n)}=\p(\J)$ for some $n$.\medskip

The name comes from the fact that the numbers ${S_{\ga,\la}\over
S_{\L^0,\la}}$  are the eigenvalues of the fusion matrix $N_\ga$, whose
entries are fusion coefficients (3.4e).

Lemma 3.3 below tells us that it suffices for most purposes to look at the
$\G$-rows and -columns of $M$. To find fusion-generators, by (3.4d) it
is natural to
look at the representation ring of $X_r$ -- see Proposition 4.1 below.
For example, if $\J=\{id.\}$, a fusion-generator is formed from the lifts
$(k-a_i^\vee)\L_0 +\L_i$ into $P_+$ of the horizontal fundamental
weights $\bar{\L}_i$, for all $i=1,\ldots,r$ (provided $k\ge {\rm max}_i
\{a_i^\vee\}$).

\medskip{{\smcap Lemma 3.3.}}\quad (a) {\it Let $\Gamma$ be a
fusion-generator for $\J$. Let $M$ be a positive invariant obeying (3.7d)
with $\J_L=\J_R=\J$, such that $(\ga,\ga)$ is $M$-monogomous for all $\ga\in
\Gamma$. Then for all $\la\in \p(\J)$, $(\la,\la)$ is $M$-monogomous.}

(b)\quad {\it Let $M$ be any positive invariant obeying (3.7d). Let
$\Gamma_L$ and $\Gamma_R$ be fusion-generators for $\J_L$ and $\J_R$,
respectively. Suppose that each $\ga\in\Gamma_L$ is right-$M$-monogomous,
and each $\ga\in\Gamma_R$ is left-$M$-monogomous.
Then there exists a map $\si:\p(\J_L)\rightarrow\p(\J_R)$ such that $(\la,
\si\la)$ is $M$-monogomous, and the induced map $\J_L\la\mapsto \J_R\si\la$
is a bijection
between the $\J_L$-orbits  in $\p(\J_L)$ and the $\J_R$-orbits in $\p(\J_R)$.
In addition,}
$$\eqalignno{S_{\la,\mu}=&\sqrt{{\|\J_R\si\la\|\,\|\J_R\si\mu\|\over
\|\J_L\la\|\,\|\J_L\mu\|}}\,S_{\si\la,\si\mu},\qquad\forall\la,\mu\in\p(\J_L)
\ ;&(3.12a)\cr
\sum_{J\in\J_L}N_{\la,\mu}^{J\nu}=&\sqrt{{\|\J_R\si\la\|\,\|\J_R\si\mu\|\,
\|\J_R\si\nu\|\over\|\J_L\la\|\,\|\J_L\mu\|\,\|\J_L\nu\|}}\,\sum_{J'\in\J_R}
N_{\si\la,\si\mu}^{J'\si\nu},\quad \forall \la,\mu,\nu\in\p(\J_L).
&(3.12b)\cr}$$
{\it Finally, let $M'$ be any other positive invariant obeying (3.7d) for the
same $\J_L,\J_R$. If for all $\ga\in\G_L$, each $(\ga,\si\ga)$ is also
$M'$-monogomous, then $M=M'$}.

\medskip{\it Proof.} \quad(a) Write $M=\oplus_i M_i$, where each $M_i$ is
irreducible (Lemma 3.2(a)). We get the equations
$$\|\J\|\,S_{\ga,\la}=(MS)_{\ga,\la}=(SM)_{\ga,\la}=\sum_\mu S_{\ga,\mu}\,
M_{\mu,\la},\quad \forall \la\in P_+\ ,\eqno(3.13)$$
and all $ \ga\in \{\L^0\}\cup \G$. In other words, writing
$v^{\ga}_i$ for the vector with components $(v^\ga_i)_\la=S_{\ga,\la}$ for
all $\la\in \I(M_i)$,  (3.13) tells us that $v^\ga_i$ is a left-eigenvector
of $M_i$ with the Perron-Frobenius eigenvalue $\|\J\|=r(M_i)$. Since $M_i$
is irreducible, this means [14] each $v^{\ga}_i$ must be a scalar multiple
of $v^{\L^0}_i$, i.e.
$${S_{\ga,\la}\over S_{\L^0,\la}}={S_{\ga,\mu}\over S_{\L^0,\mu}}\qquad
{\rm for\ all}\ \la,\mu\in\I(M_i)\ ,$$
for each $\ga\in\G$ and each submatrix $M_i$. Hence by (3.11), $(\la,\la)$
will be $M$-monogomous for each $\la\in\G^{(1)}$. Continuing recursively,
we get the desired result.

\smallskip{(b)} \quad Let $\widetilde{M}=M^T\,M/\|\J_L\|$. Then
$\widetilde{M}$ satisfies the hypotheses of
Lemma 3.3(a) with $\J=\J_R$ and $\G=\G_R$, so for each $\mu\in\p(\J_R)$,
$(\mu,\mu)$ is $\widetilde{M}$-monogomous. Thus if $\mu'\not\in\J_R\la$, then
for any $\nu\in P_+$, $M_{\nu,\mu}M_{\nu,\mu'}=0$. Similarly, for any
$\la\in\p(\J_L)$ and $\nu\in P_+$,
$M_{\la,\nu}M_{\la',\nu}=0$ whenever $\la'\not\in\J_L\la$. These
statements force the existence of the map $\si$. Equations (3.6b), (3.9) now
directly  give (3.12a), and then (3.12b) follows from (3.4e).

Finally, let $\overline{M}=M\,M'{}^{T}/\|\J_R\|$; from Lemma 3.3(a) we find
$(\la,\la)$ is $\overline{M}$-monogomous for all $\la\in\p(\J_L)$. Thus,
using (3.9), $M^T=M^T\overline{M}/\|\J_L\|=\widetilde{M}\,M'{}^T/\|\J_R\|
=M'{}^T$. \quad \QED

\medskip {\it Remark 3.1}.\quad For $A_r^{(1)}$, the case we are interested in
 in this paper, Lemma 3.1(b)(ii) forces $\J_L=\J_R$, and we can
drop the hypothesis in Lemma 3.3(b) which says that
$(\ga,\ga)$ must be $M$-monogomous for all $\ga\in\G_R$ -- this will follow
{}from the $\G_L$ hypothesis. More generally, the same thing happens
{\it whenever the number of
$\J_L$-orbits in $\p(\J_L)$ equals the number of $\J_R$-orbits in $\p(\J_R)$}.
That we can ignore $\G_R$ in this case follows from the proof of Lemma 3.3(b)
given above.

\bigskip\bigskip
{{\bf 4. The proof of the Theorem.}} \quad In this section we specialize to
the affine algebra $A_r^{(1)}$, and outline the proof of Theorem 2.1. We begin
by collecting together results particular to $A_r^{(1)}$.

Fix the algebra $A_r^{(1)}$ and the level $k$.
Let $\r=r+1$, $\k=k+\r$. The set $P_+$ of weights is given by (3.1) with all
$a_j^\vee=1$. Recall the definitions of $J\in \s_{sc}$, $C$,
and $t$, given in Section 2. $\s_{sc}$ is generated by $J$, and $\s$ by
$J$ and $C$. The Kac-Peterson matrices $S$ and $T$ are defined in
(3.3b),(3.3d). We have:
$$\eqalignno{t(J^a\la)\equiv &\,ka+t(\la)\qquad ({\rm mod\ }\r)&(4.1a)
\cr T_{J^a\la,J^a\mu}=&\,\exp[\pi i\,(-2a\,t(\la)+ka\,(\r-a))/\r]\
T_{\la,\mu}&(4.1b)\cr
S_{J^a\la,J^b\mu}=&\,\exp[2\pi i\,(b\,t(\la)+a\,t(\mu)+kab)/\r]\
S_{\la,\mu}&(4.1c)\cr
S_{C\la,\mu}=&\,S_{\la,\mu}^{*}&(4.1d)\cr
T_{C\la,C\mu}=&\,T_{\la,\mu}&(4.1e)\cr}$$
(compare (4.1b),(4.1c) with (3.5b),(3.5c) -- note that $Q_{J^d}(\la)=d\,t(\la)
/\r$). The subgroups of $\s_{sc}$ are $\J_d$, which is generated by
$J^d$ where $d$ divides $\r$. We will write $[\la]$ for the orbit of
$\la$ over $\s$, i.e.\ generated by $J$ and $C$, and $[\la]_d$
for the orbit over $\J_d$, i.e.\ generated by $J^d$. Write $\p_d=
\p(\J_d)$. Let $\F_d$
denote the fixed points of
$\J_d$. They look, schematically, like $(\mu,\ldots,\mu)$ for some
$d'$-tuple $\mu$ and some multiple $d'<\r$ of $d$.

The matrix $S$ obeys a surprising relation called {\it rank-level
duality}, related to the existence of the embedding $su(\r)\oplus su(k)
\subset su(\r k)$. Given any $\la\in P_+$, define a weight $T(\la)$ of
$A_{k-1}^{(1)}$ level $r+1$, as follows. First construct the Young diagram of
$\la$: for $1\le i\le r$, its $i$th row consists of $\sum_{j=i}^r\la_j$ boxes.
Take the transpose of this diagram; deleting all columns (if any) of length
$k$, this will be the Young diagram of some level $r+1$ $A_{k-1}^{(1)}$ weight
which we
will denote by $T(\la)$. For example, $T((k-1)\L_0+\L_\ell)=(\r-\ell)
\tilde{\L}_0+\ell\tilde{\L}_1$ for all $\ell$ -- to avoid confusion {\it we
will always put tilde's over the quantities of $A_{k-1}^{(1)}$ level} $r+1$. It
was shown in [1], by applying Laplace's determinant formula [14] to a
determinant expression for the $S$ matrix entries, that
$$S_{\la,\mu}=\sqrt{k\over\r}\, \exp[{2\pi i\over \r k}t(\la)\,
t(\mu)]\ \tilde{S}^{ *}_{T(\la),T(\mu)}\ .\eqno(4.2a)$$
Note that this map $T$ defines a bijection between $\J_1$-orbits in $P_+$
and  $\tilde{\J}_1$-orbits in $\tilde{P}_+$.

(4.2a) takes a simpler form for $\la,\mu\in
\p_1$: for $\la\in\p_1$, let $T'(\la)=\tilde{J}^{-t(\la)/\r}T
(\la)$; then the obvious calculation from (4.2a) gives
$$\eqalignno{S_{\la,\mu}=&\,\sqrt{k\over\r}\,\tilde{S}^{*}_{T'(\la),
T(\mu)}\qquad \forall \mu\in P_+&(4.2b)\cr =&\,\sqrt{k\over \r}\,
\tilde{S}^{*}_{T'(\la),T'(\mu)}\qquad \forall \mu\in\p_1\
.&(4.2c)\cr}$$
Since, in addition, $\la\in \p_1$ implies
$$\tilde{t}(T'(\la))\equiv \r\,(-t(\la)/\r)+\tilde{t}(T(\la))\equiv 0\quad
({\rm mod}\ k)\ ,$$
we see that $T'$ takes $\p_1$ to $\tilde{\p}_1$. In fact, (4.2b)
implies that $T'$ is a bijection.

Incidently, similar formulae hold for the matrix $T$. In particular,
$$\eqalignno{T_{\la,\la}\,T^{ *}_{\L^0,\L^0}=&\, \exp[{\pi i\over
\r k}t(\la)(\r k-t(\la))]\ \tilde{T}^{*}_{T(\la),T(\la)}\,\tilde{T}_{
\tilde{\L}^0,\tilde{\L}^0}\qquad \forall \la\in P_+&(4.3a)
\cr =&\, \tilde{T}^{*}_{T'(\la),T'(\la)}\,\tilde{T}_{\tilde{\L}^0,
\tilde{\L}^0}\qquad \forall \la\in\p_1\ .&(4.3b)\cr}$$

Choose any \ade\ $M$. Then by Lemma
3.1(b), $\J_L(M)$ is a subgroup of $\s_{sc}$, so equals $\J_d$ for some $d$.
Lemma 3.1(b) also says $\J_d=\J_R(M)$ -- write $\J(M)=\J_d$.
Then  by (3.6a),
$${k}'d^2\equiv 0\qquad ({\rm mod}\ 2\r)\ ,\eqno(4.4)$$
where ${k}'$ is defined in Section 2. (4.4) and (4.1b) together imply
$$(J^d\la+\rho\,|\,J^d\la+\rho)\equiv (\la+\rho\,|\,\la+\rho)\quad ({\rm mod}\
 2\r)\qquad \forall\la\in\p_d\ .$$

The first important step in the proof is to find fusion-generators (Definition
 3.3) for $\J_d$. Let
$$\eqalignno{\la^i\eqde&(k-2)\L_0+\L_i+\L_{\r-i}\qquad {\rm for}\ 1\le
i\le \r/2\ ,&(4.5a)
\cr \mu^j\eqde&(k-3)\L_0+\L_1+\L_j+\L_{r-j}\qquad {\rm for}\ 1\le j\le r/2\
,&(4.5b)\cr
\L^i\eqde& (k-1)\L_0+\L_i\qquad {\rm for}\ 0\le i<\r\ ,\qquad \L^{\r}\eqde
\L^0\ .&(4.5c)\cr}$$

\medskip{{\smcap Proposition 4.1.}}\quad (a) {\it For any divisor $d$ of
$\r$, let $\Fb_1$ denote the field
of all rational functions over $\Q$ in the Weyl characters $\overline{ch}_{
\bar{\nu}}$, for all horizontal weights $\bar{\nu}\in\cup_{k=1}^\infty
\overline{\p_d}$. Then $\Fb_d$ is generated by the
characters of the $\r$ weights $\bar{\la}^i$, $\bar{\mu}^j$, and
$\bar{\L}^{\r/d}$.}

{(b)}\quad {\it The $\r$ weights in $\Gamma_d=\{\la^i,\mu^j,\L^{\r/d}\}$
form a fusion-generator for $\J_d$. Choosing $m,\ell\in\Z$ such
that $0\le \r/d-mk=\ell<k$, the $k$ weights in $\Gamma'_d=\{\widetilde{T}'(
\tilde{\la}^i),\widetilde{T}'(\tilde{\mu}^j),J^m\widetilde{T}(\tilde{\L}^\ell)
\}$ also form a fusion-generator for $\J_d$.}\medskip

{{\it Proof.}} {(a)}\quad  $\Fb_d$ is a subfield of the field of
fractions of
the representation ring of $A_r$, the latter being isomorphic by a theorem of
Chevelley (see e.g.\ [3])
to the ring of polynomials $\Q[x_1,\ldots,x_r]$, the isomorphism sending
$$\bar{ch}_{\bar{\L}_i}\mapsto x_i\ .$$
Thus $\Fb_d$ can be thought of as the field generated over $\Q$ by the
monomials $x_1^{a_1}\cdots x_r^{a_r}$, where $t(a)\eqde\sum ja_j\equiv 0$
(mod $\r/d$). Consider the case $d=1$; the result for general $d$ follows
immediately from this.

Note that the $r$ monomials $x_1x_r,x_2x_r^2,\ldots,x_r^{\r}$ clearly
generate $\Fb_1$:
$$x_1^{a_1}\cdots x_r^{a_r}=(x_r^{\r})^{-t(a)/\r}\prod_{j=1}^r
(x_jx_r^j)^{a_j}\ .$$
Thus so do $x_i x_{\r-i}$ and $x_1 x_j x_{r-j}$: recursively,
$$x_\ell x_r^\ell=(x_{\ell-1}x_r^{\ell-1})\cdot {(x_1 x_r)\,(x_\ell
x_{\r-\ell})
\over (x_1 x_{\ell-1} x_{r-\ell+1})}\ .$$

Let $\overline{\Omega}(\bar{\mu})$ denote the dominant weights of
$\bar{L}(\bar{\mu})$. Now, $\bar{\nu}\in\overline{\Omega}(\bar{\mu})$ iff
$\bar{\mu}-\bar{\nu}$ is a sum of positive roots. From this we find
$$\eqalign{\overline{\Omega}(\bar{\la}^i)=&\,\{0,\bar{\la}^1,\ldots,
\bar{\la}^i \}\ ,\cr
\overline{\Omega}(\bar{\mu}^j)=&\,\{0,\bar{\la}^1,\ldots,\bar{\la}^{j+1},
\bar{\mu}^1,\ldots,\bar{\mu}^j\}\ .\cr}$$
Thus we can express $x_ix_{\r-i}$ and $x_1x_jx_{r-j}$ in terms of the
polynomials corresponding to $\overline{ch}_{\bar{\nu}}$, for $\bar{\mu}\in
\{\bar{\la}^i,\bar{\mu}^j\}$.\smallskip

{(b)} \quad Define $X_d=\{\la\in\p_d\,|\,S_{\L^i,\la}\ne 0\ \forall i\}$.
We begin by showing that $X_d\subseteq \G_d^{(1)}$. Indeed, choose any
$\la\in X_d$, $\mu\in\p_d$ with
$${S_{\ga,\la}\over S_{\L^0,\la}}={S_{\ga,\mu}\over S_{\L^0,\mu}}\qquad
\forall \ga\in\G_d\ .\eqno(4.6a)$$
Then by choosing $\ga=\la^i$ for various $i$ we find that also $\mu\in X_d$.
Now, for any $\alpha,\beta\in\p_d$, (3.4d) and Prop.\ 4.1(a) tell us that
$S_{\al,\be}/S_{\L^0,\be}$ can be expressed as a rational function in the
numbers $S_{\ga,\be}/S_{\L^0,\be}$, $\ga\in\G_d$; we find from the proof
of Prop.\ 4.1(a) that this expression will be well-defined (i.e.\ not of
the form 0/0) if $\be\in X_d$. What this means is that (4.6a) implies
$${S_{\al,\la}\over S_{\L^0,\la}}={S_{\al,\mu}\over S_{\L^0,\mu}}\qquad
\forall \al\in\p_d\ .\eqno(4.6b)$$
Multiplying (4.6b) by $\sum_{J\in\J_d}S^*_{\al,J\la}$ (which vanishes for
$\al\not\in \p_d$) and summing over all $\al\in P_+$, gives us
$${\|\J_d\|\over \|\J_d\la\|}\,{1\over S_{\L^0,\la}}=\sum_{J\in\J_d}{
\delta_{J\la,\mu}\over S_{\L^0,\mu}}\ ,$$
i.e.\ $\mu\in[\la]_d$. Therefore, we indeed have $X_d\subseteq\Gamma_d^{(1)}$.

Now we will show $\G_d^{(2)}=\p_d$. Choose any $\la,\mu\in\p_d$ such that
(4.6a) holds. Define the sums
$$s(\al,\be)=\,{S_{\L^0,\mu}\,S_{\L^0,\la}\over S_{\L^0,\al}\,S_{\L^0,\be}}\,
\sum_{\nu\in\p_d}S_{\nu,\al}\,S_{\nu,\be}^*\,\prod_{i=1}^r |S_{\L^i,\nu}|
\ .$$
If $s(\la,\la)\ne s(\la,\mu)$, then $\exists \nu\in\p_d$ such that $S_{\L^i,
\nu}\ne 0$ $\forall i$, and $S_{\nu,\la}^*S_{\L^0,\mu}/S_{\L^0,\la}\ne
S_{\nu,\mu}^*$ -- by the previous result this would mean $\la\in\G_d^{(2)}$.
Similarly if $s(\nu,\nu)\ne s(\la,\nu)$. Thus we may assume $s(\la,\la)
=s(\la,\nu)=s(\nu,\nu)$. But then the triangle inequality tells us
$$0=s(\la,\la)+s(\nu,\nu)-2\,s(\la,\mu)\ge \sum_{\nu\in\p_d}\bigl\{
(\sqrt{{S_{\L^0,\mu}\over S_{\L^0,\la}}}|S_{\nu,\la}|-\sqrt{{S_{\L^0,\la}
\over S_{\L^0,\mu}}}|S_{\nu,\mu}|)^2\bigr\}\prod_{i=1}^r|S_{\L^i,\nu}|\ge
0.$$
Therefore we must have $S_{\nu,\la}\,S_{\nu,\mu}^*\ge 0$ $\forall \nu\in\p_d$
(it is strictly larger for e.g.\ $\nu=\L^0$). Now from (3.4e) applied to
$N_{\L^0,\la}^{\mu'}=\delta_{\mu',\la}$, we find
$$\sum_{\mu'\in[\mu]_d}\delta_{\mu',\la}=\sum_{\mu'\in[\mu]_d}\sum_{\nu\in
\p_d}S_{\la,\nu}\,S_{\mu',\nu}^*=\|[\mu]_d\|\sum_{\nu\in\p_d}S_{\la,\nu}
S_{\mu,\nu}^*>0\ .$$
This calculation forces $\mu\in[\la]_d$, which concludes the proof that
$\G_d$ is a fusion-generator of $\J_d$.

We find that $\{\nu\in\p_d\,|\,S_{(k-i)\L_0+i\L_1,\nu}\ne 0\ \forall i\,\}
\subseteq \G_d'{}^{(1)}$, by using (4.2b) and the $\G_1$ argument. That
$\G_d'{}^{(2)}=\p_d$ now follows by the argument used for $\G_d$, with
$\L^i$ there replaced here with $(k-i)\L_0+i\L_1$.
\quad \QED \medskip

Strictly speaking, $\G_d$ requires $k\ge 3$: if $k=2$, simply drop the $\mu^j$.
We will discuss the trivial case $k=1$ at the end of this section.
Similarly, $\G_d'$ requires $r\ge 2$: if $r=1$ drop the $\widetilde{T}'(
\tilde{\mu}^j)$.

When $\r\le k$, we will choose the fusion-generator $\Gamma=\Gamma_d$;
when $\r>k>1$ we will usually choose the smaller set $\Gamma=\Gamma_d'$.
Note that $\widetilde{T}'(\tilde{\la}^i)=(k-2i)\L_0+i\L_1+i\L_r$, and
$\widetilde{T}'(\tilde{\mu}^j)=(k-2j-1)\L_0+(j-1)\L_1+\L_2+(j+1)\L_r$.

Suppose $M_{\la,\mu}\ne 0$. Then using Lemma 3.1(a),
$${\r\over d}\,S_{\L^0,\mu}=\sum_\nu M_{\L^0,\nu}\,S_{\nu,\mu}
=\sum_\nu S_{\L^0,\nu}\,M_{\nu,\mu}\ge \|\J_d\la\|\,S_{\L^0,\la}\,
M_{\la,\mu}\ . \eqno(4.7a)$$
Thus if we can show
$${\rm max}_{\ga\in\G}\{\D^{r,k}(\ga)\}<{\|\J_d \va\|\over \r/d}\,\D^{r,k}
(\va)\qquad \forall \va\in\F_d\ ,\eqno(4.7b)$$
where $\D^{r,k}(\la)$ is the q-dimension defined in (3.4b),
then that would mean $M_{\ga,\va}=M_{\va,\ga}=0$ for all fixed points $\va$ and
each $\ga\in\Gamma$. In the following section we use this idea together
with Lemmas 3.2(b) and 3.3(b) to prove:

\medskip{{\smcap Proposition 5.1$'$.}}\quad {\it For all but finitely many
pairs $(r,k)$, each $\la\in \p_d$ is right-$M$-monogomous, and neither
$\la^1$ nor $\si\la^1$ will be fixed points of $\J_d$.}\medskip

{}From the previous discussion, this is not too surprising a result
considering
that for fixed $r$, the minimum of the r.h.s.\ of (4.7b) tends to $\infty$ as
$k\rightarrow\infty$, while the l.h.s.\ tends to the dimension of some
$\overline{L}(\overline{\ga})$. The remaining $(r,k)$ are treated in section 7.

Let $\si$ denote the map in Lemma 3.3(b). Provided $\|[\la^1]_d\|=
\|[\si\la^1]_d\|$ (this will hold e.g.\ if neither $\la^1$ nor $\si\la^1$ are
fixed points), we see from (3.12a) that
$$\D^{r,k}(\la^1)=\D^{r,k}(\si\la^1)\ .\eqno(4.8)$$
Using (4.8) and (3.6a), we prove in sections 5 and 7 that:

\medskip{{\smcap Proposition 5.2$'$.}}\quad {\it Suppose all $\la\in \p_d$
are right-$M$-monogomous, and (4.8) holds. Then $\si\la^1\in[\la^1]$.}\medskip

We also know from (3.12b) that $\ga$ and $\si\ga$ have similar fusions (3.4f),
for $\ga\in\G_d$. From this we prove in section 6:

\medskip{{\smcap Proposition 6.1$'$.}}\quad {\it Suppose all $\la\in \p_d$
are right-$M$-monogomous, and (4.8) holds. Then for some $c=0,1$, we have for
all $i,j$ that
$(\la^i,C^c\la^i)$  and $(\mu^j,C^c\mu^j)$ are $M$-monogomous, and
$\si\L^{\r/d}\in[C^c\L^{\r/d}]_1$.}\medskip

Putting these all together, we are now prepared to prove:

\medskip{{\smcap Proposition 4.2.}}\quad {\it Let $M$ be an \ade\
with $\J_d=\J(M)$, and let $\Gamma$ be any fusion-generator
of $\J_d$. Suppose each $\ga\in\G$ is right-$M$-monogomous, and (4.8) holds.
Then for
some divisor $d'$ of $\r$ for which ${k}'d'$ is even, and some $c=0,1$,
we have $M=C^c\cdot I(\J_{d'})$ (see (2.3)).}\medskip

{{\it Proof.}} \quad By Remark 3.1, all $\la\in \p_d$ will be
right-$M$-monogomous. Let $c=0,1$ be as in Proposition 6.1$'$. For convenience
replace $M$ with $C^c\cdot M$. Then each $(\la^i,\la^i)$ and $(\mu^j,\mu^j)$ is
$M$-monogomous, and $(\L^{\r/d},J^m\L^{\r/d})$ is $M$-monogomous for
some $m$. By Lemma 3.3(a) and (4.4) we are done if $d=1$ (take $d'=1$), so
consider $d>1$.

Equation (3.6a), and the fact that $d$ divides $\r$, tells us
$${m\over d}+{{k}'m^2\over 2\r} \equiv 0\qquad ({\rm mod}\ 1)\ .
\eqno(4.9a)$$
Put $d'=gcd\{m,d\}$. We may assume, by adding a multiple of $d$ to $m$ if
necessary (see Lemma 3.1(a)), that $gcd\{m/d',2\r\}=1$. Then (4.9a) becomes
$${-2\r\over d}\equiv {k}'m\qquad ({\rm mod}\ 2\r/d')\ .\eqno(4.9b)$$
Consider $M'=I[\J_{d'}]$. $M'$ is a physical invariant because (2.3b)
follows from (4.9b), and the facts that $d$ divides $\r$ and $m/d'$ is odd.
Then $\J(M')=\J_d$ iff
$${2\r\over d}=gcd\{{2\r\over d'},d'{k}'\}\ .\eqno(4.9c)$$
But (4.9c) follows from (4.9b) and the fact that $d'$ divides $d$. Also,
(4.9b) tells us that $M'_{\L^{\r/d},J^m\L^{\r/d}}\ne 0$. Thus by Lemma
3.3(b), $M=M'$.\quad \QED\medskip

In particular, Propositions 4.2 and 5.1$'$  suffice to prove Theorem 2.1 for
most pairs $(r,k)$. The finitely many remaining pairs are handled in
section 7.

Incidently, Proposition 4.2 permits an immediate proof of Theorem 2.1 for
$k=1$. It suffices to note that for $k=1$: (i) $\p_d$ is generated by
$\L_{\r/d}$; (ii) $\p_d=[\L_{\r/d}]_{\r/d}$; (iii) $\p_d$ has no fixed points.
$k=1$ was first proved in [5], though in a very different way. For
convenience we will henceforth restrict attention to $k\ge 2$.

\bigskip\bigskip
{{\bf 5. Q-dimension calculations.}} \quad
The point of this section is to prove, using q-dimensions (3.4b), that for
most $(r,k)$, (4.7b) will
be satisfied, and for most $(r,k)$, (4.8) implies $\si\la^1\in[\la^1]$.
We begin by listing  some of the properties q-dimensions obey.

This section is the most technically complicated of the paper. Even so,
q-dimensions are extremely well-behaved and  amenable to
analysis. Recall their definition in (3.4b) -- the name comes from
interpreting them as ``q-deformed Weyl dimensions''.

Note from (4.1c),(4.1d) that
$$[\la]=[\mu]\quad \Longrightarrow\quad \D^{r,k}(\la)=\D^{r,k}(\mu)\ .
\eqno(5.1a)$$
Also, an immediate consequence of rank-level duality (4.2a) is that
$$\D^{r,k}(\la)=\D^{k-1,r+1}(T(\la))\ .\eqno(5.1b)$$

By $C^{r,k}$ we mean the fundamental chamber
$$C^{r,k}\eqde \{\sum_{i=0}^rx_i\L_i\,|\,x_i\in\R,\ x_i\ge 0,\ \sum_{i=0}^r
x_i=k\}\ .\eqno(5.2a)$$
Extend the domain of $\D^{r,k}$ from $P_+$ to $C^{r,k}$, using (3.4b).
Choose any real $\r$-vectors $a,b$, $b\ne 0$, such that $a+bu\in C^{r,k}$
for all $u\in [u_0,u_1]$. Then for $u_0\le u'\le u_1$, an easy calculation
gives [11]
$${d\over du}\D^{r,k}(a+bu)|_{u=u'}=0 \Rightarrow {d^2\over du^2}\D^{r,k}
(a+bu)|_{u=u'}<0\ .\eqno(5.2b)$$
This implies the important fact:
$$\forall u_0<u<u_1,\quad \D^{r,k}(a+bu)>{\rm min}\{\D^{r,k}(a+bu_0),\,
\D^{r,k}(a+bu_1)\}\ .\eqno(5.2c)$$

Elementary consequences of (5.2c) are, for all $\la\in P_+$,
$$\eqalignno{\D^{r,k}(\la)\ge&\, 1\ ,&(5.3a)\cr
\D^{r,k}(\la)=&\,1\ \quad {\rm iff}\ \quad \la\in[\L^0]\ .&(5.3b)\cr}$$

The following trigonometric identity, obtained from the factorisation of
$y^m-1$, will
be used to simplify some expressions:
$$\prod_{\ell=0}^{m-1}\sin(x+{\ell\over m}\pi)={\sin(mx)\over 2^{m-1}}\ .
\eqno(5.4)$$

Ultimately, in comparing specific q-dimensions, we will have to estimate
sizes of various products and quotients of sine's. Apart from the obvious
trigonometric identities, a useful technique is to investigate the
behaviour as $k$ or $r$ tends to $\infty$. A basic fact is that $x\,\cot x$
decreases: this means,  for $a>1$   and $0<y<x\le\pi/a$, that
$$a>\sin(ay)/\sin(y)>\sin(ax)/\sin(x)\ .\eqno(5.5a)$$
For example consider a sequence $\be^k\in C^{r,k}$ with constant projection
$\overline{\be}=\overline{\be}^k$. Then (5.5a) and the Weyl dimension formula
tells us
$$\D^{r,k}(\be^k)<\prod_{\bar{\alpha}>0}{(\overline{\be}+\bar{\rho}|\bar{
\alpha})\over (\bar{\rho}|\bar{\alpha})}={\rm dim}\,\overline{L}(\overline{\be}
)\ ,\eqno(5.5b)$$
with the l.h.s.\ converging monotonically to the r.h.s.\ as $k\rightarrow
\infty$. Another basic fact is the concavity of ln$|\sin x|$:
$$\sin(a)\,\sin(b)<\sin(a-x)\,\sin(b+x)\ ,\eqno(5.5c)$$
provided $0<b<a<\pi$ and $0<x\le(a-b)/2$.

Fix a divisor $d$ of $\r$. As in section 4, for $\r\le k$ choose the
fusion-generator $\Gamma=\Gamma_d$ (cardinality $\r$), while for $\r>k$
choose $\Gamma=\Gamma_d'$ (cardinality $k$).

\medskip{\smcap Proposition 5.1.}\quad {\it All $r\ge 1$, $k\ge 1$ satisfy
} (4.7b) {\it unless $(r,k)$ equals}

\qquad\qquad{(i)} \quad$\,\,$ $(1,k)$ {\it for} $k\in\{2,4,6,\ldots,16\}$,

\qquad \qquad{(ii)} \quad$\,$ $(2,k)$ {\it for} $k\in\{3,6,9\}$,

\qquad \qquad{(iii)} \quad$\,$ $(3,k)$ {\it for}  $k\in\{4,6,8,10\}$,

\qquad \qquad{(iv)} \quad$\,$ $(4,5)$,

\qquad \qquad{(v)} \quad$\,\,$ $(5,k)$ {\it for} $k\in\{6,8,10\}$,

\qquad \qquad{(vi)} \quad$\,\,$ (7,8),

\qquad \qquad{(vii)} \quad$\!$ $(r,6)$ {\it for} $r\in\{7,9\}$,

\qquad \qquad{(viii)} \quad$\!\!$ $(r,4)$ {\it for} $r\in\{5,7,9\}$,

\qquad \qquad{(ix)} \quad$\,\,$ $(r,3)$ {\it for} $r\in\{5,8\}$,

\qquad \qquad{(x)} \quad$\,\,\,$ $(r,2)$ {\it for} $r\in\{3,5,\ldots,15\}$.

\medskip{{\it Proof.}} \quad Choose any fixed point $\va\in \F_d$. The period
of $\va$ (with respect to $\J_d$) is
$p=d\,\|\J_d\va\|$: $\va_i=\va_j$ if $i\equiv j$ (mod $p$). Write
$$\va^p\eqde \sum_{j=0}^{\r/p-1}{kp\over \r}\,\L_{pj}\ .$$
Putting $b=J^i\va^p-J^j\va^p$ in (5.2c) and using (5.1a), we find that
$\D^{r,k}(\va)\ge \D^{r,k}(\va^p)$. Thus it suffices to consider only the
fixed points of the form $\va^p$. By (5.4) and (3.4a):
$$S_{\L^0,\va^p}=|s|\,2^{\r\,(p-1)/2}\,(\r/p)^{\r/2}\,
\prod_{j=1}^{p-1} \sin({\pi\,j\r\over p\k})^{\r-j\r/p}\ .\eqno(5.6)$$

Next, turn to the evaluation of the l.h.s.\ of (4.7b). Directly from (3.4b) we
get
$$\eqalignno{\D^{r,k}(\la^\ell)=&{\sin(\pi\,(\r-2\ell+1)/\k)\over
\sin(\pi\,(\r+1)/\k)}\,\prod_{j=1}^\ell {\sin^2(\pi\,(\r+2-j)/\k)\over
\sin^2(\pi\,j/\k)}&(5.7a)\cr
\D^{r,k}(\mu^\ell)=&{\sin(\pi\,(\r-\ell)/\k)\,\sin(\pi\,(\r+2)/\k)\,
\sin(\pi\,(\r-2\ell)/\k)\,\sin(\pi\,\ell/\k)\over \sin(\pi/\k)\,
\sin^2(\pi\,(\ell+1)/\k)\,\sin(\pi\,(\r+1)/\k)}&\cr
&\times\prod_{j=1}^\ell {\sin^2(\pi\,(\r+2-j)/\k)\over \sin^2(\pi\,j/\k)}
&(5.7b)\cr
\D^{r,k}(\L^\ell)=&\prod_{j=1}^\ell {\sin(\pi\,(\r+1-j)/\k)\over
\sin(\pi\,j/\k)}\ .&(5.7c)\cr}$$

Consider first the case $k\ge \r$.
Stirling's formula tells us $\left({m\atop \ell}\right)<2^m/\sqrt{m}$ for
all $m\ge \ell\ge 1$. Hence using (5.5b), we find
$$\D^{r,k}(\ga)<2^{2\r}\ ,\eqno(5.7d)$$
valid for all $\ga\in\G_d$, and all $\r$, $k$. By (5.5a), (5.6) and (3.4a) we
see that
$$p\,2^{-2\r}\,\D^{r,k}(\va^p)\,/\,\r\eqno(5.7e)$$
is an increasing function of $k$, for fixed $\r$ and $p$. Thus it suffices to
show (5.7e) is greater than 1 at $k=\r$.

We begin by removing the dependence of (5.7e) (at $k=\r$) on
$p$, by showing that the r.h.s.\ of (5.6) is a decreasing
function of $p=1,2,\ldots$ at $k=\r$. Indeed, using (5.4), this is
equivalent to showing
$${1\over \sqrt{2}}{\prod_{j=1}^{p-1}\sin(\pi\,j/2p)^{j/p}\over
\prod_{j=1}^p\sin(\pi\,j/(2p+2))^{j/(p+1)}}<1\ .\eqno(5.8a)$$
 Now, elementary calculus tells us $f(x)=(\sin {\pi\over 2}x
)^x$, for $0<x<1$, has a unique minimum at some point $x=x_m$ ($\approx
0.25$); $f$ decreases for $x<x_m$, and increases for $x>x_m$. Moreover,
$f(x_m)>0.78$, so $\sqrt{2}\,f(x)>1$ for all $x$. Now return
to (5.8a). Choose the $\ell\in\{0,\ldots,p-1\}$  for which ${\ell\over p}<x_m<
{\ell+1\over p}$. Because ${j\over p+1}<{j\over p}<{j+1\over p+1}$, we may
rewrite (5.8a) as
$$\bigl(\prod^\ell_{j=1}{f(j/p)\over f(j /(p+1))}\bigr)\,
\bigl({1\over \sqrt{2}\,f((\ell+1)/(p+1))}\bigr)\,
\bigl(\prod_{j=\ell+1}^{p-1}{f(j/p)\over f((j+1)/(p+1))}\bigr)
<1\ .\eqno(5.8b)$$

Thus we need only to consider $p=\r/2$ when $\r$ is even, and $p=r/2$ when
$\r$ is odd. When $\r$ is even we get from (5.6) and (5.4) that
$$\eqalignno{\D^{r,k}(\va^p)= &\,2^{\r^2/4}\prod_{{j=1
\atop j\ {\rm odd}}}^{\r-1}\sin({\pi\,j\over 2\r})^j\ge 2^{\r^2/4}\cot({\pi
\over 2\r})^{\r/2-1}\prod_{{j=1\atop j\ {\rm odd}}}^{\r-1}\sin({\pi\,j\over
2\r})^{\r/2}&\cr>&\, (.63\,\r)^{\r/2-1}2^{\r/4}&(5.8c)\cr}$$
for $\r\ge 10$, using the fact that $\cot(\pi/20)/10>0.63$ ($x\,\cot x$ is
a decreasing function). When $\r$ is odd the same argument
works, except the (5.8b) calculation is needed to replace $f({j\over
\r})$ with $f({j\over r})$ or $f({j-1\over r})$; $\r$ odd (and greater than
10) also obeys the lower bound of (5.8c) (in fact it obeys a somewhat
stronger inequality).

Collecting everything, what we have shown thus far is that (for $\r\ge 10$)
$${\|\J_d\va\|\over \|\J_d \|}\,{\D^{r,k}(\va)\over\D^{r,k}(\ga)}>
(.63\,\r)^{\r/2-1}\,2^{-7\r/4}/\r\ ,\eqno(5.8d)$$
valid for all $\va\in\F_d$, all $d<\r$ dividing $\r$, all $\ga\in\G_d$,
and all $k\ge \r$. We find that the r.h.s.\ of (5.8d) is greater than 1 for
$\r\ge 28$. To handle $\r< 28$, return to (5.7e),
find the smallest $k\ge \r$, call it $k_{r,p}$, making (5.7e) at least 1
($p<\r$ must divide $\r$, and $\r/p$ divide
$k$). In most cases $k_{r,p}=\r$; when $k_{r,p}>\r$,
tighten this estimate by going back to (5.7a), (5.7b), (5.7c), to compute
each ${p\over \r}\D^{r,k}(\va^p)/\D^{r,k}(\ga)$. We find the results given in
the statement of the proposition for $k\ge \r$.

By (5.1b), the proof for $1<k<\r$ reduces to that of $k> \r$: when $\va^p\in
P_+$, $T(\va^p)=\tilde{\va}^{kp/\r}$.\quad\QED\medskip

Proposition 5.1 and Lemmas 3.2(b), 3.3(b) give us Proposition 5.1$'$ as
stated in section 4 -- the exceptions are the pairs $(r,k)$ listed in (i)-(x).
Next we will see that in almost all cases, (4.8) forces $\si\la^1\in[\la^1]$.
\medskip

{\smcap Proposition 5.2.}\quad {\it Let ${\cal W}$ be the set of all $\nu\in
P_+$, $\nu\not\in[\la^1]$, with $\D^{r,k}(\nu)=\D^{r,k}(\la^1)$.
For any $r\ge 1$, $k\ge 2$, the only nonempty ${\cal W}$ are}:

\qquad\qquad{(a)} ${\cal W}=[\L^3]$ \quad {\it for $(r,k)=(8,3)$ and} $(8,15)$;

\qquad\qquad{(b)} ${\cal W}=[(k-3)\L_0+3\L_1]$ \quad {\it for $(r,k)=(2,9)$
and} (14,9);

\qquad\qquad{(c)} ${\cal W}=[(k-2)\L_0+2\L_2]$ \quad {\it for $(r,k)=
(3,6)$ and} (5,4);

\qquad\qquad{(d)} ${\cal W}=[\L^4]$ \quad {\it for $(r,k)=(7,4)$ and} (7,6);

\qquad\qquad{(e)} ${\cal W}=[(k-4)\L_0+4\L_1]$ \quad {\it for $(r,k)=(3,8)$
and} (5,8).

\medskip{{\it Proof.}} \quad Begin by considering, for $\r\ge 3$ and $k\ge 3$,
the quotient
$${\D^{r,k}((k-2)\L_0+\L_1+\L_2)\over \D^{r,k}(\la^1)}={\sin(\pi\,\r/\k)
\over\sin(\pi\,3/\k)}\ .\eqno(5.9a)$$
This is always $\ge 1$, with equality iff $\r=3$ or $k=3$ -- in which
case $(k-2)\L_0+\L_1+\L_2\in[\la^1]$.

Now suppose $\la\in P_+$ has at least 3 non-zero Dynkin labels $\la_i$
(so necessarily $\r,k\ge 3$).
By taking various $b=\L_i-\L_j$ in (5.2c), we find some weight $\mu$ of the
form $\mu=(k-2)\L_0+\L_m+\L_n$,
for $1\le m<n\le r$, with $\D^{r,k}(\la)\ge \D^{r,k}(\mu)$, and equality
iff $\la\in[\mu]$. Consider $T(\mu)=(\r-n)\tilde{\L}_0+(n-m)\tilde{\L}_1
+m\tilde{\L}_2$; taking $b=\tilde{\L}_0-\tilde{\L}_1$ and $b=\tilde{\L}_0
-\tilde{\L}_2$ in (5.2c) and using (5.1b), we get
$$\D^{r,k}(\mu)\ge{\rm min}\{\D^{k-1,r+1}((\r-2)\tilde{\L}_0+\tilde{\L}_1+
\tilde{\L}_2),\,\D^{k-1,r+1}(\tilde{\la}^1)\}=\D^{r,k}(\la^1)\ ,$$
using (5.9a).
Thus, $\D^{r,k}(\la)\ge\D^{r,k}(\la^1)$ for any such $\la$, with equality
iff $\la\in[\la^1]$.

On the other hand, if $\la$ has only 1 non-zero label, then $\D^{r,k}(\la)
=\D^{r,k}(\la^1)$ would require $[\la^1]=[\L^0]=[\la]$, by (5.3b).
Therefore any $\nu\in{\cal W}$ must lie in the $\s$-orbit of some
$$\nu^{ab}\eqde (k-a)\L_0+a\L_b\ .$$
We can demand $1\le a\le k/2$, and $1\le b\le \r/2$. Note that $T(\nu^{ab})=
\tilde{\nu}^{ba}$.

Consider next
$${\D^{r,k}(\nu^{22})\over \D^{r,k}(\la^1)}={\sin(\pi/\k)\,\sin^2(\pi\,\r/\k)
\over \sin^2(2\pi/\k)\,\sin(3\pi/\k)}\ .\eqno(5.9b)$$
We may suppose $k\ge 4$, $\r\ge 4$. This can be analyzed using (5.5a),
and we find that it is less than 1 iff $(\r,k)=(4,4)$, (4,5), or (5,4) and
equals 1 iff $(\r,k)=(4,6)$ and (6,4).
{}From this we deduce, in the now familiar way from (5.2c) and (5.1b), that
the only $\nu^{ab}\in{\cal W}$ with both $a>1$ and $b>1$ is $\nu^{22}$.

It thus suffices, using (5.1b), to consider $\nu^{1b}=\L^b$ for $b\le \r/2$.
$\D^{r,k}(\L^b)$ is computed in (5.7c). First note that
$$\D^{r,k}(\L^1)<\D^{r,k}(\L^2)<\cdots<\D^{r,k}(\L^{\r/2})\ .\eqno(5.10a)$$
Write ${\cal Q}_b=\D^{r,k}(\L^b)/\D^{r,k}(\la^1)$.
When $\r\le 3$ we need to consider only $\L^1$: we find that ${\cal Q}_1=1$
iff $\L^1\in[\la^1]$. Assume now that $\r\ge 4$. From (5.5c),
${\cal Q}_1<{\cal Q}_2<1$, except for $k=2$, when $\L^2\in[\la^1]$.
We find, using (5.5a), (5.5c),  that ${\cal Q}_3\ge 1$ except for
$\r=6$, 7 and 8, for all $k$, and $\r=9$ for $4\le k\le 14$. Thus from
(5.10) the only other possible $\L^b\in{\cal W}$ is $\L^4$
for $\r=8$ and any $k$, or for $\r=9$ and $4\le k\le 14$. These possibilities
are handled in the usual way.\quad\QED

\bigskip\bigskip{{\bf 6. Fusion coefficient calculations.}}
\quad
In this section we conclude the proof of Theorem 2.1 for most pairs $(r,k)$.
In particular, let $M$ be a \ade\ with $\J(M)=\J_d$.
We will use the fusion-generator $\Gamma_d$ for all $r,k$. Prop.\ 5.1$'$ says
that for each $\ga\in\Gamma_d$, there exist weights $\si\ga\in
\p_d$, $\si\ga\not\in\F_d$, such that $(\ga,\si\ga)$ is
$M$-monogomous. Proposition 5.2
tells us $\si\la^1\in[\la^1]$, for most $(r,k)$.
Our first task will be to show $\si\la^1\in[\la^1]_d$. Using this, we will
show that $\si\la^i\in[\la^i]_d$ and (replacing $M$ if necessary by $C\cdot
M$)
$\si\mu^j\in[\mu^j]_d$. From this we can obtain $\si\L^{\r/d}\in[\L^{\r/d}]_1$.
As indicated at the end of section 4, this will conclude the proof of
Theorem 2.1 for those $(r,k)$. Our main tool in this section will be {\it
fusion coefficients} (3.4e).

$\overline{L}(\bar{\la}^1)$ has dominant weights $\bar{\la}^1$ (multiplicity
1) and 0 (multiplicity $r$, as seen by e.g.\ the Weyl
dimension formula). Therefore  (3.4f) implies
$$N_{\la^1,\mu}^\nu=\left\{\matrix{n(\nu)-1&{\rm for}\ \quad\nu=\mu\cr
1&{\rm if}\ \nu=\L_i-\L_{i-1}-\L_j+\L_{j-1}+\mu\ {\rm for\ some}\ i\ne j\cr
0&{\rm otherwise}\cr}\right.\ ,\eqno(6.1)$$
where $n(\nu)$ denotes the number of Dynkin labels $\nu_i\ne 0$, and where
we put $\L_{\r}=\L_0$.
{}From (3.12b) we see that, for all $\ga,\ga',\ga''\in\G_d$,
 there is a $c(\ga,\ga',\ga'')>0$ such that
$$\sum_{J\in\J_d}N_{\ga,\ga'}^{J\ga''}=c(\ga\ga'\ga'')\sum_{J\in\J_d}N_{
\si\ga,\si\ga'}^{J\si\ga''}\eqno(6.2)$$
(each $c(\ga,\ga',\ga'')=1$ when $(r,k)$ avoids those pairs listed in
Props.\ 5.1, 5.2, but we will be more general here).

\medskip{\smcap Proposition 6.1.}\quad {\it Suppose all $\ga\in\G_d$ are
$M$-monogomous
and $\si\la^1\in[\la^1]$. Then there exists a $c\in\{0,1\}$ such that, for
all $\ga\in\G_d$, $\ga\ne\L^{\r/d}$, we have} $\si\ga\in [C^c\ga]_d.$\medskip

{\it Proof.}\quad
To begin we find from (6.1) that (defining $\la^0=\L^0$ and discarding
$\mu^0$)
$$N_{\la^1,\la^\ell}^{\nu}\ne 0\ \quad {\rm iff}\ \quad\nu\in\{\la^\ell,\,
\la^{\ell-1},\,\la^{\ell+1},\,\mu^\ell,\,\mu^{\ell-1},\,C\mu^\ell,\,
C\mu^{\ell-1},\,\la^1+\la^\ell-\L^0\}\ . \eqno(6.3)$$
Therefore $\si$ must permute the $\J_d$-orbits of these $\nu$, if it fixes
$[\la^1]_d$ and $[\la^\ell]_d$. Note also that
$$\eqalignno{(\la^i+\rho\,|\,\la^i+\rho)=&\,2i\,(\r+1-i)+(\rho\,|\,\rho)
&(6.4a)\cr
(\mu^j+\rho\,|\,\mu^j+\rho)=&\,2\r+2+2j\,(\r-j)+(\rho\,|\,\rho)&(6.4b)\cr
(\la^1+\la^h-\L^0+\rho\,|\,\la^1+\la^h-\L^0+\rho)=&\,2\r+4+2h\,(\r+1-h)
+(\rho\,|\,\rho)&(6.4c)\cr}$$
for $0\le i\le \r/2$, $0<j\le r/2$, and $0<h\le\r/2$.
The point of (6.4) is (3.6a): $\si$ must preserve norms (mod $2\k$).

Now let us show $\si\la^1\in [\la^1]_d$. (3.4e) and (4.1c) tell us that
$$N_{J^a\la,J^b\mu}^{J^{a+b}\nu}=N_{\la,\mu}^\nu,\qquad \forall \la,
\mu,\nu\in P_+\ .\eqno(6.5)$$
We know from Proposition 5.2 that $\si\la^1=J^a\la^1$, for some $a$. So,
{}from (6.1), (6.2) and (6.5) we obtain
$$0\ne\sum_{J\in\J_d}N_{\la^1,\la^1}^{J\la^1}=c(\la^1\la^1\la^1)\sum_{J\in\J_d}
N_{J^a\la^1,  J^a\la^1}^{JJ^a\la^1}=c(\la^1\la^1\la^1)\sum_{J\in\J_d}
N_{\la^1,\la^1}^{JJ^{-a}\la^1}\ .\eqno(6.6)$$
Using (6.3), (6.4) and (6.6) we get $\si\la^1\in[\la^1]_d$.

By induction on $\ell$, we find from (6.3) and (6.4) that
$\si$ must  fix each $[\la^\ell]_d$, and $[\si\mu^\ell]_d\in\{[\mu^\ell]_d,\,
[C\mu^\ell]_d\}$. Replacing $M$ if necessary with $C\cdot M$, we may suppose
that $\si$ also fixes $[\mu^1]_d$.

Now note from (6.3) that
$$\eqalignno{N_{\la^1,\mu^\ell}^{\nu}\ne 0\ \quad{\rm iff}\ \quad\bar{\nu}
\in\{&\,\overline{\mu}^\ell,\,\overline{\mu}^{\ell+1},\,
\overline{\mu}^{\ell-1},
\,\overline{\la}^{\ell},\,\overline{\la}^{\ell+1},\,\overline{\L}_2+
\overline{\L}_{\ell-1}+\overline{\L}_{r-\ell},&\cr
&\overline{\L}_2+\overline{\L}_\ell+\overline{\L}_{r-\ell-1},\,2
\overline{\L}_1+\overline{\L}_{\ell-1}+\overline{\L}_{r-\ell},
\,2\overline{\L}_1+\overline{\L}_\ell+\overline{\L}_{r-\ell-1},&\cr&
\overline{\L}_2+\overline{\L}_\ell+\overline{\L}_{r-\ell}+\overline{\L}_r,\,
\overline{\la}^1+\overline{\la}^\ell,\,\overline{\la}^1+\overline{\la
}^{\ell+1},\,\overline{\la}^1+\overline{\mu}^\ell\,\}\ .&(6.7)}$$
Label these weights consecutively $\nu^1=\mu^\ell,\ldots,\nu^{13}=\la^1+
\mu^\ell-\L^0$. When we must make $\ell$ explicit, we will write these as
$\nu^a(\ell)$. Assume inductively that $\si\mu^1=\mu^1,\ldots,\si\mu^\ell
=\mu^\ell$, and suppose for contradiction that $\si\mu^{\ell+1}=C\mu^{\ell+1}$.
We are interested here in $k\ge 3$, $\r\ge 5$, $1\le\ell\le {r\over 2}-1$.
The question is, when can $C\mu^{\ell+1}\in[\nu^a(\ell)]_d$ for some $a\ge 6$?

Solving this is straightforward once one realizes that we can ignore
$\nu^6(1)$, $\nu^6(2)$, $\nu^7(1)$, $\nu^8(1)$, $\nu^{10}(1)$ -- in these cases
$\nu^a(\ell)$ equals one of $\la^{\ell+1},\mu^\ell,\mu^{\ell+1},$ or $\la^1
+\la^{\ell+1}-\L^0$. This means that the $\bar{\L}_i$ for each $\bar{\nu}$
in (6.7) are written in non-decreasing order of indices, making comparison
with $J^mC\mu^{\ell+1}$ easy. Note also that we can ignore $\nu^{11}$ and
$\nu^{12}$ for all $\ell$, as $C\mu^{\ell+1}$ lies in their $\J_d$-orbit
iff $\mu^{\ell+1}$ does.
Also, for $k\ge 7$ we find that $(C\mu^{\ell+1})_0=k-3>\nu_i^a$ for $i>0$.
So the only possibilities are $3\le k\le 6$.

One solution is $k=4$, $\r=12$, and $\ell=3$, where we have $C\mu^{\ell+1}
=J^{-3}\nu^{10}$. But $\si\mu^4\ne C\mu^4$ because both $d=1$ and $d=3$
violate (4.4). Similarly, when $k=\r=6$ and $\ell=1$, we have $C\mu^2=J^{-1}
\nu^{13}$, but $d=1$ violates (4.4).

The only other solution is $k=3$, $\ell=\r/3$: we find $C\mu^{\ell+1}=J^{\r/3}
\nu^6$. However, if $\si\mu^{\ell+1}=C\mu^{\ell+1}$ in this case, then by
(6.7) and (3.12b) we would need to have $\mu^{\ell+1}\in[\nu^a(\ell-1)]_d$
for some $a,d$, and this does not happen here.

Thus in all cases $\si\mu^{\ell+1}\in[\mu^{\ell+1}]_d$, and our proposition
is proved.\quad\QED\medskip

{}From (3.12a), we find that Proposition 6.1 tells us
$${S_{\ga,\la}\over S_{\L^0,\la}}={S_{\ga,\si\la}\over S_{\L^0,\si\la}}
\quad \forall\ga\in\G_d,\ \ga\ne \L^{\r/d},\ \forall\la\in\p_d\ .\eqno(6.8a)$$
But these $\ga$ form a fusion-generator for $\J_1$, so (6.8a) implies
$\si\la\in[\la]_1$ for all $\la\in \p_1$, by the recursive argument in the
proof of Lemma 3.3(a). By the argument given in the proof
of Prop.\ 4.1(b), we also see that (6.8a) implies
$${S_{\la,\mu}\over S_{\L^0,\mu}}={S_{\la,\si\mu}\over S_{\L^0,\si\mu}}
\quad \forall\la\in\p_1\ ,\eqno(6.8b)$$
and all $\mu\in X_d$. These two observations, together with the fact (see
the proof of Prop.\ 4.1(b)) that $\la\in X_d^{(1)}=\p_d$, imply $\si\la\in
[\la]_d$ for all $\la\in\p_1$. Hence by (3.12a), (6.8b) holds in addition
for all $\mu\in\p_d$. Multiplying (6.8b) by $\sum_{J\in\J_1}S_{\la,J\mu}^*$
(which vanishes for $\la\not\in\p_1$) and summing over all $\la\in P_+$,
we find $\si\mu\in[\mu]_1$. In particular this holds for $\mu\in\L^{\r/d}$,
which gives us Proposition 6.1$'$ as stated in section 4.

\bigskip\bigskip{{\bf 7. Anomolous ranks and levels.}}\quad
In this section we conclude the proof of Theorem 2.1 by addressing the few
pairs $(r,k)$ which slipped through
the previous arguments. In the process we will explicitly construct the
${\cal E}_7$-type exceptionals listed in section 2. Recall
Proposition 4.2.

\medskip{{\it 7.1.}}\quad In this subsection we use norm arguments and
tighten the q-dimension arguments to discard almost all remaining pairs
$(r,k)$.

Consider first Proposition 5.2. It is trivial to verify that (3.6a) is
violated by  the choice $\la=\la^1$, $\mu=\nu$ for any $\nu$ and $(r,k)$
listed in (a)-(e). This gives us Proposition 5.2$'$ stated in section 4.

Next we turn to (i)-(x) in Proposition 5.1. Again (3.6a) is a severe
constraint, as is the requirement (see (4.7a)) that
$$M_{\la,\mu}\ne 0\Rightarrow {\|\J_d\|\over \|\J_d\la\|}\ge {\D^{r,k}(\la)
\over \D^{r,k}(\mu)}\ge{\|\J_d\mu\|\over\|\J_d\|}\ .\eqno(7.1a)$$
Moreover, if in addition $\|\J_d\|/\|\J_d\mu\|<(\D^{r,k}(\mu)+1)/\D^{r,k}(\la)$
 then
$$M_{\la,\nu}=\delta_{[\nu]_d,[\mu]_d}\qquad \forall \nu\in P_+\ .
\eqno(7.1b)$$

Once again we will choose the fusion-generator $\G=\G_d$ when $\r\le k$, and
$\G=\G_d'$ when $\r>k>1$.
Using (3.6a), (7.1a), (4.4) and Proposition 5.1, we find that
$M_{\ga,\va}=0$ for all $\ga\in\G$ and $\va\in \F_d$, except
possibly for:

\item{(a)$'$} $(r,k,d)=(1,4,1)$, $\ga=\va=\la^1$;

\item{} $(r,k,d)=(1,16,1)$, $\ga=\la^1$ and $\va=(8,8)$;

\item{(b)$'$} $(r,k,d)=(2,3,1)$, $\ga=\va=\la^1$;

\item{} $(r,k,d)=(2,9,1)$, $\ga=\la^1$ and $\va=(3,3,3)$;

\item{(c)$'$} $(r,k,d)=(3,2,2)$, $\ga=\va=\la^1$;

\item{} $(r,k,d)=(3,4,2)$, $\ga=\va=\la^2$;

\item{} $(r,k,d)=(3,8,d)$ for $d=1$ or 2, $\ga=\la^1$ and $\va\in[(4,0,4,0)
]_1$, and (if $d=1$) $\ga=\mu^1$ and $\va=(2,2,2,2)$;

\item{(d)$'$}   $(r,k,d)=(4,5,1)$, $\ga=\la^1$ and $\va=(1,1,1,1,1)$;

\item{(e)$'$} $(r,k,d)=(5,6,2)$, $\ga=\mu^1$, $\va=(2,0,2,0,2,0)$;

\item{(f)$'$} $(r,k,d)=(7,4,d)$ for $d=2$ or 4, $\ga=\la^1$ and $\va\in[
2\L_0+2\L_4]_2$, and (if $d=2$) $\ga=\widetilde{T}'\tilde{\mu}^1$ and
$\va=\L_1+\L_3+\L_5+\L_7$;

\item{(g)$'$} $(r,k,d)=(8,3,3)$, $\ga=\la^1$ and $\va=\L_0+\L_3+\L_6$;

\item{(h)$'$} $(r,k,d)=(15,2,4)$, $\ga=\la^1$ and $\va=\L_0+\L_8$;

\item{} $(r,k,d)=(15,2,8)$, $\ga=\la^1$ and $\va\in[\L_0+\L_8]_4$.

\medskip{{\it 7.2.}}\quad
In this subsection we find all the ${\cal E}_7$-type exceptionals
for $A_1^{(1)}$. This was first done in [4].
By the previous subsection it suffices to consider $(k,d)=(4,1)$ and (16,1).

First consider $k=4$. The problem here is that $\la^1$ is a fixed point.
However, $\p_1=[\L^0]_1\cup[\la^1]_1$, so $(\la^1,\la^1)$ must be
$M$-monogomous. Thus there is no exceptional here.

Next, consider $k=16$. Write $\va=(8,8)$. Suppose $M$ is exceptional. Then
by Proposition 4.2, both $M_{\la^1,\va}\ne 0$ and $M_{\va,\la^1}\ne 0$.
By (7.1b) we find that in fact
$$M_{\la^1,\nu}=M_{\nu,\la^1}=\delta_{\nu,\va},\qquad \forall \nu\in P_+
\ .\eqno(7.2a)$$
The remaining entries of $M$ are fixed by Lemma 3.1, (3.9) and (3.6a), except
for
$M_{\va,\va}$. $M_{\va,\va}=1$ is forced by Lemma 3.2(a): we must have the
eigenvalue
$$r_m\eqde r(\left(\matrix{m&1&1\cr 1&0&0\cr 1&0&0\cr}\right) )\eqno(7.2b)$$
equal to $\r/d=2$; by [14], $r_m$ is a strictly increasing function of
$m\ge 0$, and an easy calculation gives $r_1=2$.

\medskip{{\it 7.3.}}\quad In this subsection we find all the ${\cal E}_7$-type
exceptionals for $A_2^{(1)}$. This was first done in [10].
By subsection 7.1 it suffices to consider $(k,d)=(3,1)$ and (9,1). $k=3$
is handled by the identical argument used on $k=4$ in subsection 7.2.

So consider $k=9$, and assume $M$ is exceptional. Write $\va=(3,3,3)$. Exactly
as before we must have (7.2a)
satisfied. There are precisely 7 $\J_1$-orbits in $\p_1$; (3.6a) and
Lemma 3.2(b) force $(\la,\la)$ to be $M$-monogomous for each $\la\in
\p_1$  except for $\la\in[\la^1]\cup [\va]\cup[(0,3,6)]$;
either $(\la,\la)$ or $(\la,C\la)$ will be $M$-monogomous for each
$\la\in[(0,3,6)]$. The value of $M_{\va,\va}$ again is fixed by the (7.2b)
argument. Note that
$((0,3,6),\,(0,3,6))$ being $M$-monogomous implies $M={\cal E}^{(2,9)}$, while
$((0,3,6),\,C(0,3,6))$ being $M$-monogomous implies $M=C\cdot
{\cal E}^{(2,9)}$.

\medskip{{\it 7.4.}}\quad In this subsection we find all the ${\cal E}_7$-type
exceptionals for $A_3^{(1)}$.
By the previous subsection it suffices to consider $(k,d)=(2,2)$, (4,2), (8,2)
and (8,1). For both $k=2$ and $k=4$, the problem is that
$\la^{k/2}$ is  a fixed point. However in both cases (3.6a) forces $(
\la^{k/2},\la^{k/2})$ to be $M$-monogomous, and so there are no exceptionals.

Consider next $k=8$ and $d=2$, and suppose $M$ is exceptional. Then
$M_{\la^1,\va}\ne 0$ for some $\va\in[(4,0,4,0)]_1$.
By (7.1b) we again have the analogue of (7.2a) satisfied.
We find, by (3.6a) and (7.1a), that $(2\rho,2\rho)$ must be $M$-monogomous.
We can show (see e.g.\ (3.4a)) that
$$S_{2\rho,\la^1}=-S_{2\rho,\va}\ne 0\ .\eqno(7.3)$$
However, this violates (3.6b) evaluated at $(\la^1,2\rho)$. Thus
$M$ cannot be exceptional.

Finally, consider $k=8$ and $d=1$, and let $M$ be exceptional. $\p_1$
consists of exactly 12
$\J_1$-orbits; (3.6a) and (7.1a) tell us  that the only way to have both
$M_{\la,\mu}\ne 0$ and $\la\not\in[\mu]_1$ is if either $\la,\mu\in[\mu^1]_1
\cup[C\mu^1]_1\cup[2\rho]_1$ or $\la,\mu\in[\la^1]_1\cup[\va]_1$, where
$\va=(4,0,4,0)$. Because $M$ commutes with $C=S^2$, we see that
$M_{\mu^1,2\rho}=M_{C\mu^1,2\rho}$. By (7.1b) we find that $M_{\mu^1,2\rho}
\ne 0$ iff
$$M_{\mu,\la}=M_{\la,\mu}=\delta_{\la,2\rho}, \qquad \forall \mu\in[\mu^1],
\quad \forall \la\in P_+\ .\eqno(7.4)$$
When $M_{\mu^1,2\rho}\ne 0$, we fix the value $M_{2\rho,2\rho}=2$ by the
usual argument as in (7.2b). Similarly, if $M_{\la^1,\va}\ne 0$, then we
know all $[\la^1]\cup[\va]$-rows and -columns of $M$.

The only remaining question is whether $M_{\mu^1,2\rho}\ne 0$ iff
$M_{\la^1,\va}
\ne 0$. That this is so follows immediately from (7.3), and the fact that
(2.7a) will satisfy (3.6b). Hence $M={\cal E}^{(3,8)}$.

\medskip{{\it 7.5.}}\quad In this subsection we find all the
${\cal E}_7$-type exceptionals for $A_4^{(1)}$ at $k=5$, $d=1$.
There are precisely 6 $\J_1$-orbits in $\p_1$, and exactly one fixed
point: $\va=\rho$. The argument is exactly as for $(r,k,d)=(1,16,1)$, and
we get that $M={\cal E}^{(4,5)}$.

\medskip{{\it 7.6.}}\quad In this subsection we show there are no
${\cal E}_7$-type exceptionals for $A_5^{(1)}$ at $k=6$, $d=2$. Put
$\va=(2,0,2,0,2,0)$. By (e)$'$ and Prop.\ 4.2, any exceptional $M$
would have both $M_{\mu^1,\va}$ and $M_{\va,\mu^1}$ nonzero. Put $\la=4\L_0
+\L_1+\L_2$; by the usual arguments we find that in fact $(\la,J^aC^b\la)$ is
$M$-monogomous for some $a,b\in\{0,1\}$. Replacing $M$ with
$I(\J_3)^a \cdot C^b\cdot M$ ($I(\J_3)$ is an invertible matrix), we may
assume $(\la,\la)$
is $M$-monogomous. From (3.6b) we get $3\,S_{\la,\mu^1}=S_{\la,\va}$.
However its l.h.s.\ is
non-real, while the r.h.s.\ is real. So such an $M$ cannot exist.

\medskip{{\it 7.7.}}\quad In this subsection we find all the
${\cal E}_7$-type exceptionals for $A_7^{(1)}$.
By the previous subsection it suffices to consider $(k,d)=(4,4)$ and (4,2).
The case $k=d=4$ is handled exactly as $(r,k,d)=(3,8,2)$ was (the analogue
of (7.3) here is simply its rank-level dual -- see (4.2c)).

When $(k,d)=(4,2)$, there are precisely 24 $J$-orbits, including four
orbits of fixed points. Equations (7.1a) and (4.6a) require $(J\L^0,
J\L^0)$ to be $M$-monogomous, so by Lemma 3.1(a) $M_{\la,\mu}=M_{J\la,
J\mu}$ for all $\la,\mu$. The remainder of the argument is exactly as
for $(r,k,d)=(3,8,1)$.

\medskip{{\it 7.8.}}\quad In this subsection we find all the
${\cal E}_7$-type exceptionals for $A_8^{(1)}$, when $k=d=3$. There are 21
$\J$-orbits, and three fixed points. We find from (7.1a) and (4.6a) that,
replacing $M$ with $C\cdot M$ if necessary,
$(J\L^0,J\L^0)$ must be $M$-monogomous, so by Lemma 3.1(a) we know
that $M_{\la,\mu}=M_{J\la,J\mu}$ for all $\la,\mu\in P_+$, and
$$M_{\la,\mu}\ne 0\Rightarrow t(\la)\equiv t(\mu)\qquad ({\rm mod}\ 9)\ .
\eqno(7.5)$$
Applying (7.5), (7.1a) and (4.6a) to the remaining weights in $\p_3$,
and using the familiar arguments, we can fix all values of $M$, except for
determining whether $(\mu^1,\mu^1)$ or $(\mu^1,C\mu^1)$
is $M$-monogomous. However the former must hold, because otherwise $M-
{\cal E}^{(8,3)}+I(\J_3)$ would violate Proposition 4.2.

\medskip{{\it 7.9.}}\quad In this subsection we find all the
${\cal E}_7$-type exceptionals for $A_{15}^{(1)}$, when $k=2$ and $d=8$ or 4.
Start with $d=8$. There are
40 $\J_4$-orbits, with 8 containing fixed points. But everything simplifies,
as we once again find (conjugating if necessary) that
$M_{\la,\mu}=M_{J\la,J\mu}$ and
$$M_{\la,\mu}\ne 0\Rightarrow t(\la)\equiv t(\mu)\qquad ({\rm
mod}\ 16)\eqno(7.6)$$
for all $\la,\mu\in P_+$ by the usual arguments. We can now quickly force
$M={\cal E}^{(15,2)}$.

The case $d=4$ is completely analogous to $(r,k,d)=(1,16,1)$.

\bigskip\bigskip{{\bf 8. Conclusion.}} \quad
The problem of classifying all {\it physical invariants} (see Definition 3.1)
for a given nontwisted affine algebra $X_r^{(1)}$ and level $k$ is a key
step in the classification problem for RCFTs.
Though this problem remains open, progress is being made, and with this
paper the end could be in sight, at least for the algebras $A_r^{(1)}$.
In particular, there is a natural division of the problem into two
pieces, based on the structure of the physical invariant about
the distinguished weight $k\L_0$. One subproblem is to classify all physical
invariants $M$ -- they are called {\it \ade s} -- whose
$k\L_0$-row and column reflect the symmetries of the Dynkin diagram of
$X_r^{(1)}$ (see (1.3b) for a more precise statement).
Almost all physical invariants are expected to be
of this type, including what seem to be the most elusive {\it exceptional}
physical invariants. In this paper we develop an approach for achieving the
classification of \ade s for any algebra $X_r^{(1)}$, and explicitly solve it
for the algebra $A_r^{(1)}$ at all levels $k$.
We find that most of these physical invariants are built up in a natural
way from the symmetries of the Coxeter-Dynkin diagram, but some
are not -- the so-called ${\cal E}_7$-{\it type exceptionals}. These
exceptional invariants are surprisingly rare.

The second subproblem, which remains completely open apart from some minor
special cases
(most notably $A_1^{(1)}$ [4] and $A_2^{(1)}$ [10]), is to find those
$X_r^{(1)}$ and $k$ for which a physical invariant can have `irregular'
values of $M_{k\L_0,\la},M_{\la,k\L_0}$.
In the language of conformal field theory, this is the problem of finding
all possible exceptional {\it chiral extensions} for the given affine algebra
$X_r^{(1)}$ at level $k$. These also seem to be quite rare: for $A_r^{(1)}$ at
level $k$, these occur at $(r,k)=
(r,r-1)$, $(r,r+1)$, and $(r,r+3)$, and probably only finitely many other
pairs ($r,k)$. The solution to these two subproblems would quickly imply the
classification of all physical invariants.

A natural follow-up to this paper will be to extend the results here to
the remaining affine algebras. $A^{(1)}_r$ is special because its
Coxeter-Dynkin
diagram is so symmetrical. In some ways this makes the classification of
\ade s  more difficult (e.g.\ we are inundated
with possibilities at most steps of the solution), but in other ways it
makes things much simpler (e.g.\ q-dimensions are easier to handle).
In any case, it can be expected that the classification for the remaining
$g=X^{(1)}_r$ should follow the method used here [13].

Another follow-up is the classification of all physical invariants for
$A_r^{(1)}$ at levels 2 and 3 (until now only level 1 is known). This
follows from Theorem 2.1 in this paper, and the work in [10], and will be
reported elsewhere.

Of course it would be preferable for a ``uniform proof'' of the classification
of \ade s for all $X_r^{(1)}$ and $k$. At
present the closest we can come are the lemmas in section 3. In fact it
is far from clear that even a ``uniform'' {\it list} of \ade s is
 possible -- the problem of course are the ${\cal E}_7$-type exceptionals.

Another disappointing feature of the proof given here is that it is
 computer-assisted (although in a minor way): some q-dimensions (3.4b)
were computed for sections 5 and 7, and an $S$ matrix element
was needed in section 7.6. Round-off error is not a problem in
these calculations (2 decimal places suffice). Nevertheless it would
be nice if these computations could be replaced by more conceptual arguments.
 Again, the existence of the ${\cal E}_7$-type exceptionals makes this more
difficult.

This paper hints that the classification of all physical invariants, at
least for simple $X_r$, may not be far away.
In the process of solving this problem we are being
forced to investigate properties of the Kac-Peterson matrices in some
detail, and it can be hoped that our analysis will also be of use to other
problems involving the modular behaviour of the affine characters.

\bigskip{{\it Acknowledgments.}}\quad It is a pleasure to thank the
hospitality of IHES, the Concordia University math department, and MPIM,
the places where
this paper was written. Conversations with Philippe Ruelle and Mark Walton
are also appreciated.\bigskip

\bigskip \centerline{{\smcap References}} \medskip

\item{[1]} {\smcap D.\ Altschuler, M.\ Bauer, and C.\ Itzykson}, {\it The
branching rules of conformal embeddings,} Commun.\ Math.\ Phys.\ {\bf 132}
(1990), 349-364.

\item{[2]} {\smcap D.\ Bernard}, {\it String characters from Kac-Moody
automorphisms,} Nucl.\ Phys.\ {\bf B288} (1987), 628-648.

\item{[3]} {\smcap N.\ Bourbaki}, {\it Groupes et Alg\`ebres de Lie,
Chapitres IV-VI,} Hermann, Paris, 1968.

\item{[4]} {\smcap A.\ Cappelli, C.\ Itzykson and J.-B.\ Zuber}, {\it The
A-D-E classification of $A_1^{(1)}$ and minimal conformal field theories,}
Commun.\ Math.\ Phys.\ {\bf 113} (1987), 1-26.

\item{[5]} {\smcap P.\ Degiovanni}, {\it $Z/NZ$ conformal field theories},
Commun.\ Math.\ Phys.\ {\bf 127} (1990), 71-99.

\item{[6]} {\smcap C.\ Dong, H.\ Li and G.\ Mason}, {\it Simple currents
and extensions of vertex operator algebras}, preprint (q-alg/9504008).

\item{[7]} {\smcap G.\ Felder, K.\ Gawedzki and A.\ Kupiainen}, {\it
Spectra of Wess-Zumino-Witten
models with arbitrary simple groups,} Commun.\ Math.\ Phys.\ {\bf 117} (1988),
127-158.

\item{[8]} {\smcap A.\ Font}, {\it Automorphism fixed points and
exceptional modular invariants}, Mod.\ Phys.\ Lett.\ {\bf A6} (1991),
3265-3272.

\item{[9]} {\smcap T.\ Gannon}, {\it Towards a classification of
su(2)$\,\oplus\cdots\oplus\,$su(2) modular invariant partition functions,}
J.\ Math.\ Phys.\ {\bf 36} (1995), 675-706.

\item{[10]} {\smcap T.\ Gannon}, {\it The
classification of SU(3) modular invariants revisited}, to appear in
Annales de l'I.\ H.\ P., Phys.\ Th\'eor.

\item{[11]} {\smcap T.\ Gannon}, {\it Symmetries of the Kac-Peterson modular
matrices of affine algebras,} to appear in Invent.\ Math.\ (q-alg/9502004).

\item{[12]} {\smcap T.\ Gannon, Ph.\ Ruelle, M.\ A.\ Walton}, {\it
Automorphism modular invariants of current algebras,} to appear in
Commun.\ Math.\ Phys.\ (hep-th/9503141).

\item{[13]} {\smcap T.\ Gannon}, {\it The automorphisms of simple current
chiral extensions} (work in progress).

\item{[14]} {\smcap F.\ R.\ Gantmacher}, {\it The Theory of Matrices,}
Chesea Publishing Co., New York, 1990.

\item{[15]} {\smcap K.\ Gawedzki}, {\it Conformal Field Theory}, in
Ast\'erisque {\bf 177-178} (1988/89), 95-126.

\item{[16]} {\smcap V.\ G.\ Kac,} {\it Infinite Dimensional Lie algebras,}
3rd edition, Cambridge University Press, Cambridge, 1990.

\item{[17]} {\smcap V.\ G.\ Kac and D.\  Peterson}, {\it Infinite-dimensional
Lie algebras, theta functions and modular forms,} Adv.\ Math.\ {\bf 53}
(1984), 125-264.

\item{[18]} {\smcap V.\ G.\ Kac and M.\ Wakimoto}, {\it Modular and conformal
constraints in representation theory of affine algebras,} Adv.\ Math.\
{\bf 70} (1988), 156-236.

\item{[19]} {\smcap M.\ Kreuzer and A.\ N.\ Schellekens,} {\it Simple currents
versus orbifolds with discrete  torsion -- a complete classification,} Nucl.\
Phys.\ {\bf B411} (1994), 97-121.

\item{[20]} {\smcap G.\ Moore and N.\ Seiberg}, {\it Naturality in conformal
field theory,} Nucl.\ Phys.\ {\bf B313} (1989), 16-40.

\item{[21]} {\smcap A.\ Ocneanu}, {\it Topological QFT and subfactors},
talks at Fields Institute, April 26-30, 1995.

\item{[22]} {\smcap Ph.\ Ruelle, E.\ Thiran and J.\  Weyers}, {\it
Implications of an arithmetic symmetry
of the commutant for modular invariants,} Nucl.\ Phys.\ {\bf B402} (1993),
693-708.

\item{[23]} {\smcap A.\ N.\ Schellekens and S.\ Yankielowicz}, {\it Modular
invariants from simple currents. An explicit proof,} Phys.\ Lett.\ {\bf B227}
(1989), 387-391.

\item{[24]} {\smcap A.\ N.\ Schellekens and S.\ Yankielowicz}, {\it Field
identification fixed points in the coset construction}, Nucl.\ Phys.\
{\bf B334} (1990), 67-102.

\item{[25]} {\smcap G.\ Segal}, {\it The definition of conformal field theory},
Oxford University preprint, 1987;

\item{} {\it Geometric aspects of quantum field theories}, in: {\it Proc.\ of
the ICM, Kyoto}, Springer-Verlag, Hong Kong, 1991.

\item{[26]} {\smcap P.\ Slodowy}, {\it A new A-D-E classification}, Bayreuther
Math.\ Schr.\ {\bf 33} (1990), 197-213.

\item{[27]} {\smcap M.\ A.\ Walton}, {\it Algorithm for WZW fusion rules: a
proof,} Phys.\ Lett.\ {\bf B241} (1990), 365-368.

\item{[28]} {\smcap J.-B.\ Zuber}, {\it Graphs and reflection groups},
to appear in Commun.\ Math.\ Phys.\ (hep-th/9507057).

\end